\documentclass[trsc]{informs4}
\makeatletter
\def\theARTICLETOPRIGHT{}
\makeatother
\usepackage{eqndefns-left} 
\RequirePackage{tgtermes}
\RequirePackage{newtxtext}
\RequirePackage{newtxmath}
\RequirePackage{bm}
\RequirePackage{endnotes}

\OneAndAHalfSpacedXII 
\usepackage{multirow}
\usepackage{algorithm}
\usepackage{algpseudocode}
\usepackage{tikz}
\usepackage{microtype}
\usepackage{graphicx}
\usepackage{subcaption}     
\usepackage{booktabs}       
\usepackage{array}
\usepackage{tabularx}
\usepackage{longtable}

\usepackage{amsmath}
\usepackage{amssymb}
\usepackage{mathtools}
\usepackage{enumitem}
\usepackage{fontawesome5}
\usepackage{float}
\usepackage{hyperref}
\usepackage[capitalize,noabbrev]{cleveref}

\usepackage{natbib}
 \bibpunct[, ]{(}{)}{,}{a}{}{,}%
 \def\bibfont{\small}%
\EquationsNumberedThrough    

\TheoremsNumberedThrough     
\ECRepeatTheorems  %

\MANUSCRIPTNO{TRSC-0001-2024.00}

\begin{document}



\RUNAUTHOR{Abdellaoui et al.}

\RUNTITLE{Learning implicit feasibility constraints}

\TITLE{Learning Implicit Feasibility Constraints for Real-World Routing and Scheduling: Application to Log Transportation}

\ARTICLEAUTHORS{%
\AUTHOR{Abdelhakim Abdellaoui *}
\AFF{Department of Mathematics and Industrial Engineering,
Polytechnique Montreal, \EMAIL{abdelhakim.abdellaoui@polymtl.ca}}
\AFF{Natural Ressources Canada, CanmetEnergy-Varennes}
\AUTHOR{Ayoub Boufous}
\AFF{Department of Mathematics and Industrial Engineering,
Polytechnique Montreal}
\AUTHOR{Issmail El Hallaoui}
\AFF{Department of Mathematics and Industrial Engineering,
Polytechnique Montreal}
\AUTHOR{Loubna Benabbou}
\AFF{Department of management science, Université du Québec à Rimouski (UQAR)}
\AUTHOR{Francois Aubé}
\AFF{Natural Ressources Canada, CanmetEnergy-Varennes}
\AUTHOR{Mouloud Amazouz}
\AFF{Natural Ressources Canada, CanmetEnergy-Varennes}
} 

\ABSTRACT{%

Real-world vehicle routing and scheduling problems involve complex operational rules and feasibility constraints typically formulated as mixed-integer linear programs (MILP). However, optimization tools are built around a fixed set of hard-coded constraints, while in practice this set evolves as new rules or preferences emerge, seasonally or permanently. Updating it requires modeling and operations research skills that planners rarely have, so generated plans are routinely adjusted by hand based on practical knowledge. Building on recent work that uses machine learning to recover such hidden constraints, we propose a data-driven constraint-learning approach that trains three complementary predictors, a Graph Neural Network (GNN), a decision tree, and a linear regression, on historical execution data from a log-truck routing and scheduling problem ($\mathcal{LTRSP}$), and embeds each inside a MILP through linearized constraints. We further introduce a stacking mechanism that combines all three within a single augmented optimization problem (AOP), letting the solver endogenously select the most reliable predictor for each decision. On real-world industrial data, each predictor already improves feasibility, but the stacked embedding consistently achieves the lowest objective degradation: it (i)~satisfies the operational rules on unseen instances with smaller degradation than any single-model variant, (ii)~picks the most appropriate predictor per decision without prior knowledge of the rule's nature, and (iii)~reduces daily manual adjustment effort while remaining tractable for daily use. Beyond this application, the framework enables optimization tools that adapt to evolving practice without recurrent manual remodeling.

}%




\maketitle

\section{Introduction and related work}

Optimization techniques are widely used in routing and scheduling applications to assist planners in making operational decisions \citep{laporte2009fifty}. In practice, these systems commonly rely on MILP formulations that encode a predefined set of constraints intended to represent the operational environment \citep{toth2014vehicle}. While such models are effective in capturing formal business rules and resource limitations, real-world vehicle routing and scheduling problems rarely obey only the constraints explicitly written in the planning model.

In operational settings, planners routinely follow additional unwritten rules, context-dependent conventions, temporal habits, resource-specific compatibilities, and tacit feasibility restrictions. These implicit operational constraints arise from evolving practices, local knowledge, and external factors such as weather or infrastructure conditions \citep{pawlak2017automatic}. As a result, the true feasible region of the routing system is often significantly smaller than that assumed by the optimization model, leading to solutions that are mathematically feasible but operationally rejected or significantly modified before execution \citep{hewitt2020data}.

MILP solutions generated by state-of-the-art solvers are frequently modified manually before execution to reflect practical considerations \citep{bayani2024learning}. These manual adjustments, however, are typically performed locally and under time pressure, without access to the global optimization context. Consequently, they may introduce inefficiencies and inconsistencies, and can lead to solutions that deviate from optimality and vary in quality across planners and planning horizons. The resulting discrepancies between optimized and executed plans contain valuable information about the true operational feasibility of the system, and have motivated a growing body of work on learning from historical decisions to improve optimization models.

In general, the discrepancies between the optimal plan and the executed plan arise from two fundamentally different sources. The first consists of deterministic operational rules that planners apply systematically and recurrently before execution, reflecting operational knowledge not encoded in the model; these are the implicit constraints that our framework aims to learn. The second consists of stochastic disruptions occurring during execution, such as demand fluctuations, road closures, or equipment failures, which also force plan revisions but are random, instance-specific, and non-recurrent. Disruptions of this kind are not implicit constraints and are addressed by reoptimization and recovery methods; they fall outside the scope of this work.

The problem of recovering hidden information from observed decisions has a long history in operations research under the framework of inverse optimization, which aims at inferring the parameters of an optimization problem from solutions believed to be optimal for an unobserved decision-maker \citep{ahuja2001inverse,heuberger2004inverse}. Recent extensions have significantly broadened this paradigm: data-driven inverse optimization handles multiple noisy observations and provides statistical guarantees \citep{aswani2018inverse, esfahani2018datadriven}, while contextual and learning-based variants infer parameters as functions of side information \citep{bertsimas2018datadriven, sadana2025survey}. Closely related streams include constraint acquisition, which automatically infers the constraints of a model from labeled examples of valid and invalid solutions \citep{bessiere2017constraint, beldiceanu2016modelseeker, kumar2022learning}, and the broader predict-then-optimize literature, which integrates prediction and optimization in a single decision pipeline \citep{elmachtoub2022smart, vanderschueren2022predict}.

Within this landscape, two recent contributions are particularly relevant to our work and define the methodological direction in which we situate ourselves. The work of \cite{hewitt2020data} introduced a data-driven framework for customizing optimization models from discrepancies between optimized and executed plans, using linear regression to infer missing operational rules. \cite{bayani2024learning} extended this idea by relying on decision trees, which offer richer expressiveness for non-linear and logical patterns and improve interpretability for practitioners. Both works demonstrate that useful operational rules can be recovered from historical optimization-execution pairs and re-injected into the MILP through linear encodings of the learned predictor. While effective for feature-based and low-dimensional patterns, these approaches share three important limitations. First, they commit to a single predictor family, which restricts their ability to capture the heterogeneous nature of implicit operational rules encountered in practice. Second, they do not explicitly exploit the relational and graph-structured nature of decision making problems such as routing and scheduling. Third, they have been validated only on synthetic or small-scale benchmark instances, leaving open the question of whether such data-driven constraint learning approaches can scale to and remain effective in real-world industrial routing and scheduling contexts.

Routing and scheduling problems, however, are inherently combinatorial and graph-structured, with feasibility often governed by interactions between nodes, arcs, resources, and time \citep{toth2014vehicle}. In particular, the $\mathcal{LTRSP}$ studied by \cite{abdellaoui2025decomposition} is formulated over a non-acyclic time-space network, in which constraints couple routing and scheduling decisions across multiple arcs, vehicles, and time intervals. In practice, the implicit constraints that emerge from operational behavior can be of very different natures: some are relational, involving dependencies between routing entities across the graph; others are conditional, following if-then logic tied to specific operational contexts; and others are approximately linear, reflecting proportional trends in resource usage or demand patterns. No single predictor family is equally suited to all these forms. Graph Neural Networks (GNNs) are naturally adapted to capture relational dependencies through message passing over graph topology \citep{scarselli2008graph, wu2020comprehensive, cappart2023combinatorial}, decision trees excel at encoding conditional if-then decision rules through axis-aligned partitions of the feature space, and linear regressors are well suited to capture global linear trends. The exact embedding of such predictors within a MILP has itself become an active research area, with formal encodings developed for trained neural networks \citep{fischetti2018deep, anderson2020strong, grimstad2019relu} and unified through software frameworks such as OMLT \citep{ceccon2022omlt}.

In this work, we study the problem of learning and integrating implicit operational constraints in large-scale $\mathcal{LTRSP}$. Building on the methodological lineage of \cite{hewitt2020data} and \cite{bayani2024learning}, but departing from their reliance on a single predictor family, we develop a framework in which three complementary predictors, a GNN, a decision tree, and a linear regression, are each trained on historical executed routing and scheduling plans to learn feasibility patterns. Each trained model is then embedded exactly into a MILP formulation through a linear encoding, leveraging recent advances in MILP-compatible representations of machine learning models \citep{fischetti2018deep, anderson2020strong, ceccon2022omlt}. Following the no-free-lunch principle, we further introduce a stacking mechanism that combines all three embedded predictors within a single augmented optimization problem (AOP) for $\mathcal{LTRSP}$. The resulting feasible region reflects both classical routing and scheduling constraints and the learned operational rules.

Using large-scale routing data from an industrial forest partner, we show that while each individual predictor already improves plan feasibility, the stacked embedding consistently incurs the lowest objective degradation while satisfying the operational rules, and substantially reduces the need for daily manual intervention by planners. Because the nature of an implicit rule, whether relational, conditional, or linear, is not known in advance, the stacking mechanism removes the need to commit to a single predictor family and provides a robust choice across the heterogeneous rules encountered in practice. More broadly, our framework advances a paradigm for optimization-based decision support in which implicit operational rules are continuously learned from adjusted plans and embedded into the model. By reducing reliance on manual constraint engineering and repeated customization, the approach enables more autonomous and resilient optimization systems, in which learned feasibility structures complement rather than replace classical mathematical modeling. Our contributions are as follows:
\begin{enumerate}
    \item We extend the data-driven constraint learning paradigm of \cite{hewitt2020data} and \cite{bayani2024learning} by introducing three complementary predictors, a GNN, a decision tree, and a linear regression, each suited to a different nature of implicit constraint encountered in large-scale routing and scheduling operations.
    \item We embed each trained predictor exactly inside a MILP through linearized constraints, preserving the exactness and optimality guarantees of the solver, and in particular provide a graph-aware embedding of a message-passing GNN that exploits the routing network topology.
    \item We introduce a confidence-based stacking mechanism that combines all embedded predictors into a single augmented routing and scheduling problem, allowing the MILP to endogenously select, for each decision, the most appropriate predictor without prior knowledge of the rule's nature. We provide a formal characterization showing that this construction is guaranteed not to underperform the best single predictor, at both the relaxed and integer levels.
    \item We evaluate the framework on a large-scale industrial dataset from a Canadian forest partner, showing that the stacked embedding achieves the lowest objective degradation across a representative set of operational rules while remaining computationally tractable for daily operational use.
\end{enumerate}
\section{Background and methodology}
Let $\mathcal{G} = (\mathcal{N}, A)$ be a directed graph representing a routing and scheduling network \cite{abdellaoui2025decomposition}, where $\mathcal{N}$ denotes the set of nodes and $A \subseteq \mathcal{N} \times \mathcal{N}$ is the set of admissible directed arcs between nodes. Each node $n \in \mathcal{N}$ is defined as a tuple $n = (\ell, \tau)$ pairing a location $\ell$ (e.g., a forest block, mill, or a homebase) with a time interval $\tau$. Each arc $a \in A$ is associated with a nonnegative cost $c_a$, typically representing distance, travel time, or transportation cost. A routing and scheduling problem consists in selecting a subset of arcs forming feasible routes that minimize the total routing and scheduling cost, subject to a set of structural, operational, and logical constraints . Over the past decades, a wide range of solution approaches have been proposed for routing and scheduling problems in the operations research literature. These approaches can be broadly classified into metaheuristic algorithms \cite{ghotb2024optimization}, exact optimization methods \cite{melchiori2022mathematical}, and hybrid approaches \cite{audy2023planning}. Exact methods typically rely on MILP formulations and aim at guaranteeing optimality \cite{abdellaoui2025decomposition}. In what follows, we present the general framework of our approach. The proposed framework is called \textbf{CLIF}, in reference to Constraint Learning for Implicit Feasibility.

\subsection{\textbf{CLIF} general framework}
\label{sec:CLIF_framework}

For a MILP formulation, such as the one proposed for $\mathcal{LTRSP}$ in \cite{abdellaoui2025decomposition} and detailed in Appendix~\ref{LTRSP}, the decision is made over keys. A key $\alpha = (a, v, t) \in \mathcal{A}$ associates an arc $a \in A$ with a vehicle $v$ and a day $t$, and $\mathcal{A}$ denotes the set of all keys. We introduce a binary routing and scheduling decision variable for each key:
\[
x_{\alpha} \;=\;
\begin{cases}
1, & \text{if key } \alpha \text{ is selected }, \\
0, & \text{otherwise},
\end{cases}
\qquad \forall\, \alpha \in \mathcal{A}.
\]
We collect these variables into the vector $\mathbf{x} = (x_{\alpha})_{\alpha \in \mathcal{A}} \in \{0,1\}^{|\mathcal{A}|}$, and we denote the solution space as $\mathcal{X} = \{0,1\}^{|\mathcal{A}|}$. Using this notation, the routing and scheduling model $\mathcal{LTRSP}$ can be formulated as follows:
\[
\begin{aligned}
\min_{\mathbf{x}} \quad & \mathbf{c}^{\top}\mathbf{x} \\
\text{s.t.} \quad & \mathbf{A}\,\mathbf{x} \;\ge\; \mathbf{b}, \\
& \mathbf{x} \in \mathcal{X} .
\end{aligned}
\]
where $\mathbf{A}$ is the constraint matrix, $\mathbf{c}$ is the cost vector, and $\mathbf{b}$ is the right-hand-side vector. We denote this generic MILP by $\mathcal{LTRSP}_G$.

In practice, this $\mathcal{LTRSP}_G$ is embedded in decision optimization software and deployed for use by operational planners. In other words, this software presents the optimal solution of $\mathcal{LTRSP}_G$ above to the planners. Let $\mathbf{x}^{\star} \in \mathcal{X}$ denote this optimal routing and scheduling plan. As mentionned before, this plan is rarely executed as-is. Instead, prior to execution, planners apply manual modifications to $\mathbf{x}^{\star}$, yielding an executed plan $\hat{\mathbf{x}} \in \mathcal{X}$ that reflects additional operational considerations not captured in the original constraint matrix $\mathbf{A}$. The systematic discrepancies between $\mathbf{x}^{\star}$ and $\hat{\mathbf{x}}$ reveal the presence of implicit constraints that restrict the set of keys effectively usable in operations. We formalize this through an implicit feasibility function $f^{\star}: \mathcal{A} \rightarrow \{0,1\}$, where $f^{\star}(\alpha)=1$ if key $\alpha$ is observed in executed plans and $f^{\star}(\alpha)=0$ otherwise. This implicit feasibility function captures systematic deviations between optimal planned solutions and their corresponding adjusted counterparts, reflecting latent operational rules that are not explicitly encoded in the original constraint set.

Our objective is to learn an approximation $\hat{f}_{\theta}(\alpha) \approx f^{\star}(\alpha)$ from historical data, capturing how operational feasibility emerges from the interaction between arcs, nodes, and their surrounding routing context, as revealed by repeated planning and execution cycles. We model each key-level feasibility prediction as $\hat{x}_{\alpha} = \hat{f}_{\theta}^{m}(\mathcal{G}, X,\alpha)$, where $\hat{f}_{\theta}^m$ is a trained predictor that estimates the expected operational value of key $\alpha$ using machine learning model $m$ given the routing and scheduling graph $\mathcal{G}$ and its feature matrix $X$. To integrate the learned feasibility function into the original $\mathcal{LTRSP}_G$, we introduce deviation variables $\Delta_{\alpha} \ge 0$ that measure the discrepancy between the decision variables $x_{\alpha}$ and their predicted feasibility $\hat{x}_{\alpha}$:
\[
\Delta_{\alpha} \ge x_{\alpha} - \hat{x}_{\alpha}, \qquad
\Delta_{\alpha} \ge \hat{x}_{\alpha} - x_{\alpha}, \qquad \forall \alpha \in \mathcal{A}.
\]
The resulting augmented optimization problem, denoted $\mathcal{LTRSP}_{\text{AUG}}^m$, such that $m$ is the predictor used to approximate $\hat{f}_{\theta}$. This model balances the original objective with learned operational feasibility by minimizing a weighted combination of routing costs and prediction deviations:
\[
\begin{aligned}
\min_{\mathbf{x},\,\hat{\mathbf{x}},\,\boldsymbol{\Delta}} \quad
& \sum_{\alpha \in \mathcal{A}} c_{\alpha} x_{\alpha} + \lambda \sum_{\alpha \in \mathcal{A}} \Delta_{\alpha} \\
\text{s.t.} \quad
& \mathbf{A}\,\mathbf{x} \;\ge\; \mathbf{b}, \\
& x_{\alpha} \in \{0,1\}, \qquad \forall \alpha \in \mathcal{A}, \\
& \hat{x}_{\alpha} = \hat{f}_{\theta}^{m}(\mathcal{G}, X,\alpha), \qquad \forall \alpha \in \mathcal{A}, \\
& \Delta_{\alpha} \ge x_{\alpha} - \hat{x}_{\alpha}, \qquad \forall \alpha \in \mathcal{A}, \\
& \Delta_{\alpha} \ge \hat{x}_{\alpha} - x_{\alpha}, \qquad \forall \alpha \in \mathcal{A}, \\
& \Delta_{\alpha} \ge 0, \qquad \forall \alpha \in \mathcal{A}.
\end{aligned}
\]
where $\lambda \geq 0$ is a penalty coefficient that controls the trade-off between the original objective and the alignment with learned feasibility patterns: a higher $\lambda$ forces the solver to produce plans closer to what the predictor considers operationally feasible, while $\lambda = 0$ recovers the original formulation. The learned feasibility function effectively reshapes the feasible region explored by the solver while preserving the tractability of the formulation.

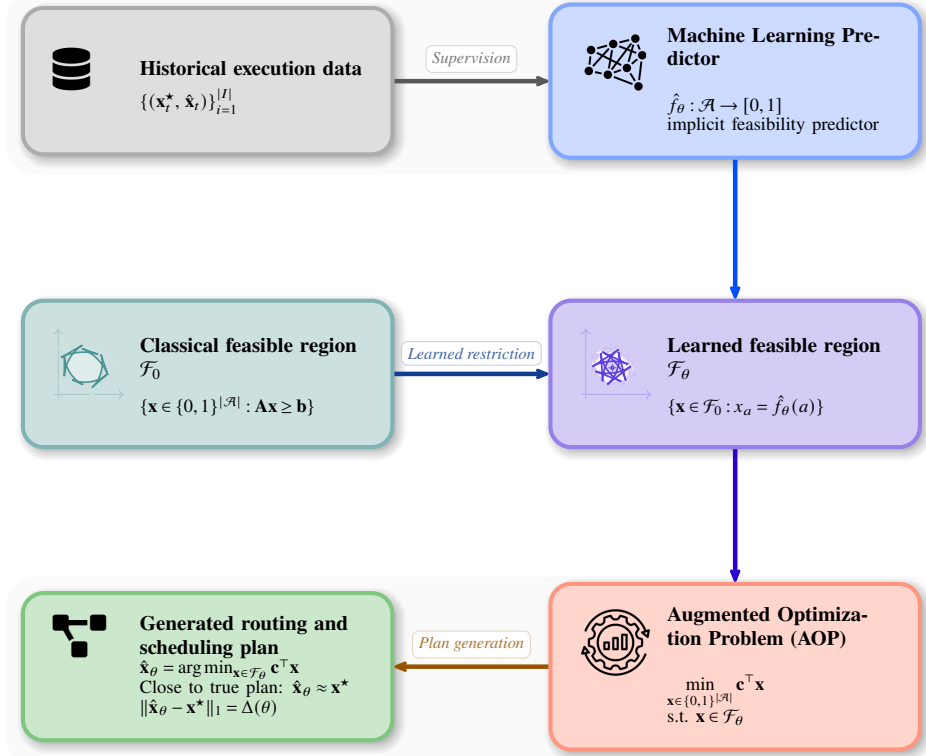
\begin{figure*}[h]
\centering
\usetikzlibrary{
  arrows.meta,
  positioning,
  shapes,
  calc,
  fit,
  backgrounds,
  decorations.pathreplacing,
  shadows.blur
}

\colorlet{cdata}{gray!70!black}
\colorlet{cgnn}{blue!70!cyan}
\colorlet{czero}{teal!80!black}
\colorlet{ctheta}{blue!65!violet}
\colorlet{caop}{red!65!orange}
\colorlet{cplan}{green!55!black}

\resizebox{0.75\textwidth}{!}{%
\begin{tikzpicture}[
  font=\small,
  mainarrow/.style={
    -{Stealth[length=9pt, width=7pt, round]},
    line width=2.2pt,
    color=#1
  },
  box/.style={
    rectangle,
    rounded corners=10pt,
    draw=#1!50,
    fill=#1!6,
    line width=2pt,
    inner sep=10pt,
    align=center,
    minimum width=6.0cm,
    minimum height=2.5cm,
    text width=5.6cm,
    blur shadow={shadow blur steps=8,
                 shadow xshift=2pt, shadow yshift=-2pt,
                 shadow blur radius=5pt,
                 fill=#1!20}
  },
  lbl/.style={
    font=\scriptsize\itshape,
    fill=white,
    inner sep=3pt,
    rounded corners=3pt,
    text=#1!80,
    draw=#1!25,
    line width=0.5pt
  }
]


\node[box=cdata] (Data) at (0,0) {
  \begin{tabular}{@{}m{1.2cm}m{3.8cm}@{}}
    {\LARGE\faDatabase} &
    \textbf{Historical execution data}\\[-2pt]
    & \footnotesize $\{(\mathbf{x}_t^\star,\,\hat{\mathbf{x}}_t)\}_{i=1}^{|I|}$
  \end{tabular}
};

\node[box=cgnn] (GNN) at (9,0) {
  \begin{tabular}{@{}m{1.2cm}m{3.8cm}@{}}
    \raisebox{-0.6cm}{\includegraphics[height=1.2cm]{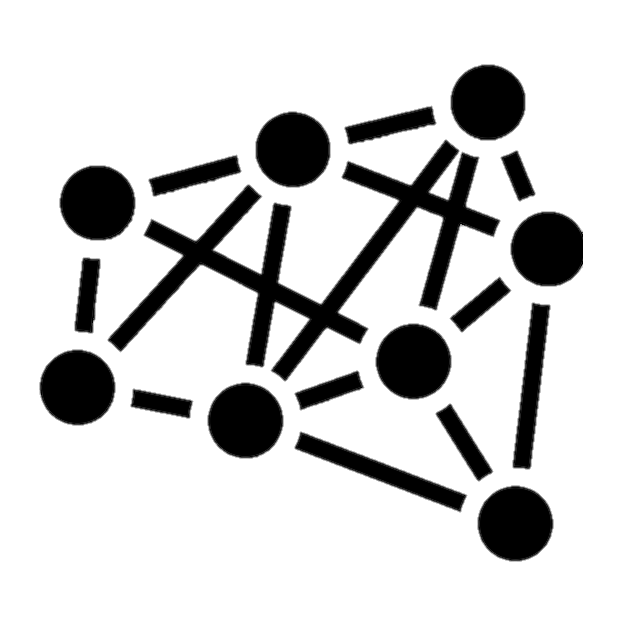}} &
    \textbf{Machine Learning Predictor}\\[-2pt]
    & \footnotesize $\hat{f}_\theta : \mathcal{A} \to [0,1]$\\[-2pt]
    & \footnotesize implicit feasibility predictor
  \end{tabular}
};

\draw[mainarrow=cdata] (Data.east) -- (GNN.west)
  node[lbl=cdata, midway, above=4pt] {Supervision};


\node[box=czero] (F0) at (0,-5) {
  \begin{tabular}{@{}m{1.2cm}m{3.8cm}@{}}
    \raisebox{-0.3cm}{%
      \begin{tikzpicture}[scale=0.7]
        \draw[czero!35, line width=0.4pt, ->] (-0.1,0) -- (1.6,0);
        \draw[czero!35, line width=0.4pt, ->] (0,-0.1) -- (0,1.6);
        \draw[czero!18, line width=0.3pt] (-0.05,1.15) -- (1.55,0.7);
        \draw[czero!18, line width=0.3pt] (0.95,1.55) -- (1.45,0.0);
        \draw[czero!18, line width=0.3pt] (-0.05,0.3) -- (0.95,-0.05);
        \fill[czero!12] (0.15,0.45) -- (0.15,0.95) -- (0.5,1.15)
          -- (1.05,0.95) -- (1.2,0.5) -- (0.8,0.15) -- cycle;
        \draw[czero!70, line width=0.8pt] (0.15,0.45) -- (0.15,0.95)
          -- (0.5,1.15) -- (1.05,0.95) -- (1.2,0.5) -- (0.8,0.15) -- cycle;
      \end{tikzpicture}%
    } &
    \textbf{Classical feasible region} $\mathcal{F}_0$\\[-2pt]
    & \footnotesize $\{\mathbf{x}\in\{0,1\}^{|\mathcal{A}|} : \mathbf{A}\mathbf{x} \ge \mathbf{b}\}$
  \end{tabular}
};

\node[box=ctheta] (Ftheta) at (9,-5) {
  \begin{tabular}{@{}m{1.2cm}m{3.8cm}@{}}
    \raisebox{-0.3cm}{%
      \begin{tikzpicture}[scale=0.7]
        \draw[ctheta!35, line width=0.4pt, ->] (-0.1,0) -- (1.6,0);
        \draw[ctheta!35, line width=0.4pt, ->] (0,-0.1) -- (0,1.6);
        \fill[ctheta!4] (0.15,0.45) -- (0.15,0.95) -- (0.5,1.15)
          -- (1.05,0.95) -- (1.2,0.5) -- (0.8,0.15) -- cycle;
        \draw[ctheta!25, line width=0.5pt, dashed] (0.15,0.45) -- (0.15,0.95)
          -- (0.5,1.15) -- (1.05,0.95) -- (1.2,0.5) -- (0.8,0.15) -- cycle;
        \fill[ctheta!18] (0.4,0.5) -- (0.38,0.82) -- (0.58,0.95)
          -- (0.85,0.82) -- (0.92,0.55) -- (0.7,0.38) -- cycle;
        \draw[ctheta!75, line width=0.8pt] (0.4,0.5) -- (0.38,0.82)
          -- (0.58,0.95) -- (0.85,0.82) -- (0.92,0.55) -- (0.7,0.38) -- cycle;
        \fill[ctheta!90] (0.65,0.65) circle(1.5pt);
        \draw[ctheta!60, line width=0.3pt] (0.65,0.55) -- (0.65,0.75);
        \draw[ctheta!60, line width=0.3pt] (0.55,0.65) -- (0.75,0.65);
      \end{tikzpicture}%
    } &
    \textbf{Learned feasible region} $\mathcal{F}_\theta$\\[-2pt]
    & \footnotesize $\{\mathbf{x}\in\mathcal{F}_0 : x_{a} = \hat{f}_\theta(a)\}$
  \end{tabular}
};

\draw[mainarrow=czero!60!ctheta] (F0.east) -- (Ftheta.west)
  node[lbl=czero!60!ctheta, midway, above=4pt] {Learned restriction};


\node[box=cplan] (Plan) at (0,-10) {
  \begin{tabular}{@{}m{1.2cm}m{3.8cm}@{}}
    {\LARGE\faProjectDiagram} &
    \textbf{Generated routing and scheduling plan}\\[-2pt]
    & \footnotesize $\hat{\mathbf{x}}_\theta = \arg\min_{\mathbf{x}\in\mathcal{F}_\theta} \mathbf{c}^\top \mathbf{x}$\\[-2pt]
    & \footnotesize Close to true plan:\enspace $\hat{\mathbf{x}}_\theta \approx \mathbf{x}^\star$\\[-2pt]
    & \footnotesize $\|\hat{\mathbf{x}}_\theta - \mathbf{x}^\star\|_1 =\Delta(\theta)$
  \end{tabular}
};

\node[box=caop] (AOP) at (9,-10) {
  \begin{tabular}{@{}m{1.2cm}m{3.8cm}@{}}
    \raisebox{-0.6cm}{\includegraphics[height=1.2cm]{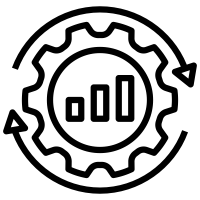}} &
    \textbf{Augmented Optimization Problem (AOP)}\\[-2pt]
    & \footnotesize $\displaystyle\min_{\mathbf{x}\in\{0,1\}^{|\mathcal{A}|}} \mathbf{c}^\top \mathbf{x}$\\[-2pt]
    & \footnotesize s.t.\ $\mathbf{x} \in \mathcal{F}_\theta$
  \end{tabular}
};

\draw[mainarrow=cgnn] (GNN.south) -- (Ftheta.north);

\draw[mainarrow=ctheta] (Ftheta.south) -- (AOP.north);

\draw[mainarrow=caop!60!cplan] (AOP.west) -- (Plan.east)
  node[lbl=caop!60!cplan, midway, above=4pt] {Plan generation};

\begin{pgfonlayer}{background}
  \fill[gray!4, rounded corners=14pt]
    ($(Data.north west)+(-0.25,0.25)$) rectangle
    ($(GNN.south east)+(0.25,-0.25)$);
  \fill[gray!4, rounded corners=14pt]
    ($(Plan.north west)+(-0.25,0.25)$) rectangle
    ($(AOP.south east)+(0.25,-0.25)$);
\end{pgfonlayer}

\end{tikzpicture}
}
\vspace{-0.5em}
\caption{
A predictor learns an implicit feasibility map from historical
execution data and is embedded exactly into a mixed-integer program.
The resulting augmented optimization problem induces a refined feasible
region $\mathcal{F}_\theta$ that converges statistically to the true
operational region $\mathcal{F}^\star$.
The generated routing and scheduling plan $\hat{\mathbf{x}}_\theta$ closely approximates the
true optimal plan $\mathbf{x}^\star$.
}
\label{feasability-map1}
\end{figure*}

This formulation is general and independent of the choice of predictor. In this work, $\hat{f}_{\theta}^{m}$ can be instantiated using either a linear regressor (LR), a decision tree (DT), a Graph Neural Network (GNN), or a stacking combination of all three (STK). The key requirement is that the predictor can be encoded as a set of linear and integer constraints within the $\mathcal{LTRSP}$. The complete \textbf{CLIF} framework is illustrated in Figure~\ref{feasability-map1}. We will subsequently detail the learning component of this framework.
The
$\mathcal{LTRSP}$-compatible embedding of LR and DT is presented in Appendix~\ref{AOP-trees-lineair}
and follows the approach of \cite{hewitt2020data,bayani2024learning}. In what
follows, we present the GNN-based embedding and the stacking
mechanism.

\subsection{GNN embedding}
\label{sec:gnn_embedding}

Among the three predictors, the GNN is suited to capture relational feasibility patterns, as it reasons over the graph topology of the routing and scheduling network through message passing. We first describe the GNN architecture and the notation used throughout, and then show how its forward pass is encoded exactly as a set of linear and integer constraints within the augmented optimization problem.

\subsubsection{GNN architecture}
\label{subsec:gnn}

\begin{figure*}[t]
  \centering
  \includegraphics[width=\textwidth]{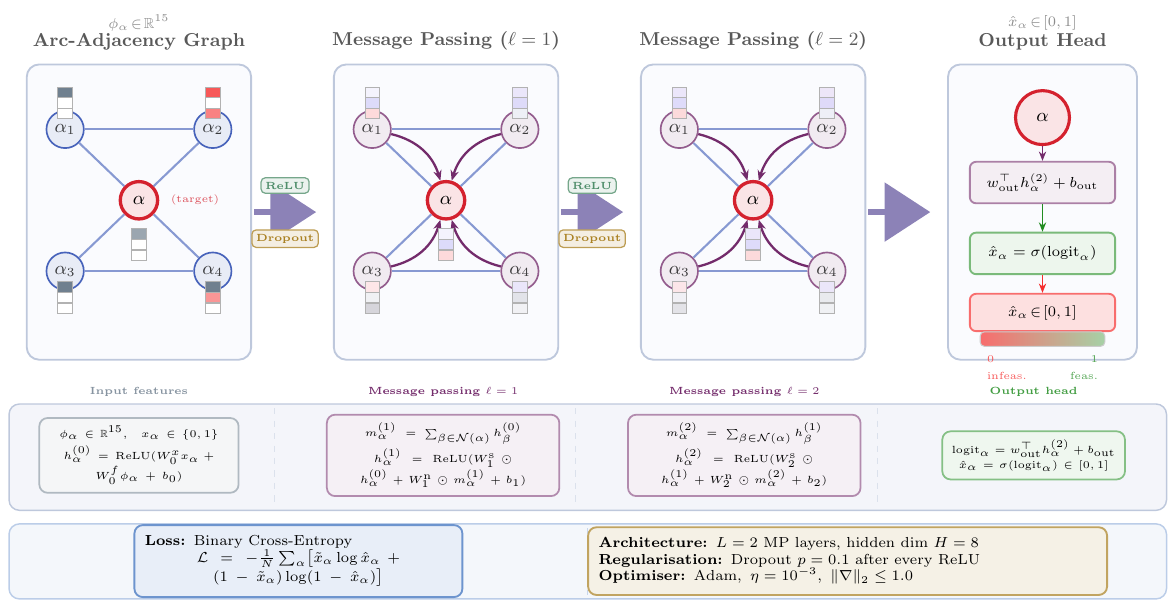}
  \caption{GNN architecture of \textbf{CLIF}. Each panel shows the arc-adjacency graph at a successive processing stage. Coloured embedding cells below each node visualise the evolution of node representations $h^{(l)}_i$ across layers. The formula summarises the update equations at every stage.}
\label{fig:gnn_architecture}
\end{figure*}

The GNN applies a single shared set of parameters to every key, so that a model trained on historical instances can be applied to new instances of different size and topology, while neighborhood aggregation lets it exploit the relational structure of the corresponding graph. We train a message-passing GNN constructed over an arc-adjacency graph, where the model parameters are learned through a supervised learning problem on graph-structured data using historical pairs of optimized and executed plans. The complete GNN architecture as well as the used notation are illustrated in Figure~\ref{fig:gnn_architecture}.

Explicitly, a node of the GNN corresponds to a key $\alpha = (a, v, t)$, where $a = (u, w)$ is an arc identified by its two endpoints, $v$ is a vehicle, and $t \in \{0,\dots,T-1\}$ is a time index. Two keys $\alpha_i = (a_i, v, t)$ and $\alpha_j = (a_j, v, t)$ are connected if their arcs $a_i$ and $a_j$ share at least one endpoint, that is, if $a_i$ and $a_j$ are adjacent in the graph. Let $N(\alpha)$ denote the neighbor set of key $\alpha$. From an optimization point of view, each node is associated with a decision variable $x_{\alpha} \in \{0,1\}$ describing whether the key is present in the original solution. From a machine learning perspective, each key is associated with a fixed feature vector $\phi_{\alpha} \in \mathbb{R}^{R}$, with $R = 15$, as illustrated in Table~\ref{feature-tabtab}. The features include normalized indices for time and vehicle, one-hot encodings of the types of the arc endpoints (home base $HB$, forest $F$, or mill $M$), normalized endpoint identifiers, bucketed endpoint identifiers, normalized departure and arrival time intervals, and a constant bias term.

\begin{table}[H]
\centering
\small
\begin{tabular}{cll}
\toprule
Index & Feature & Description \\
\midrule
$f_1$     & $t/T_{\text{size}}$               & Normalized time index \\
$f_2$     & $\mathrm{idx}(v)/V_{\text{size}}$ & Normalized vehicle index \\
$f_{3:5}$ & $\phi(u)$                         & One-hot type of origin node $u$ \\
$f_{6:8}$ & $\phi(w)$                         & One-hot type of destination node $w$ \\
$f_9$     & $\eta(u)$                         & Normalized identifier of $u$ \\
$f_{10}$  & $\eta(w)$                         & Normalized identifier of $w$ \\
$f_{11}$  & $\beta(u)$                        & Bucketed identifier of $u$ \\
$f_{12}$  & $\beta(w)$                        & Bucketed identifier of $w$ \\
$f_{13}$  & $i_{\text{from}}/I_{\max}$        & Normalized departure time interval \\
$f_{14}$  & $i_{\text{to}}/I_{\max}$          & Normalized arrival time interval \\
$f_{15}$  & $1$                               & Bias term \\
\bottomrule
\end{tabular}
\caption{
15-dimensional feature vector for key $\alpha = (a, v, t)$ with $a = (u, w)$.
The mapping $\phi(\cdot)$ denotes a 3-dimensional one-hot encoding of the node type
($\mathrm{HB},\mathrm{F},\mathrm{M}$).
The normalized identifier $\eta(\cdot)$ rescales node indices within each type.
The bucketed identifier $\beta(\cdot)$ corresponds to a modular encoding
$\mathrm{id}(n) \bmod B$ (with $B=5$) normalized to $[0,1]$.
The departure and arrival intervals $i_{\text{from}}$ and $i_{\text{to}}$
are the time slots at which the arc begins and ends, normalized by $I_{\max}=20$.
}
\label{feature-tabtab}
\end{table}

More precisely, let $x^\star_\alpha \in \{0,1\}$ denote the value of key $\alpha$ in the optimal plan, used as a scalar input feature of the GNN. The model first computes a hidden representation $h^{(0)}_{\alpha} \in \mathbb{R}^{H}$ for each node $\alpha$ using an input embedding:
\[
h^{(0)}_{\alpha}=\mathrm{ReLU}\!\left( W^{\alpha}_{\mathrm{0}}\,x^\star_\alpha 
+ W^{f}_{\mathrm{0}} \phi_\alpha + b_{\mathrm{0}} \right),
\]
where $W^{\alpha}_{\mathrm{0}}\in\mathbb{R}^{H}$ multiplies the scalar $x^\star_\alpha$, $W^{f}_{\mathrm{0}}\in\mathbb{R}^{H\times R}$, and $b_{\mathrm{0}}\in\mathbb{R}^{H}$. The network then performs $L$ message-passing layers. At layer $l$, neighbor aggregation is defined as
\[
m^{(l)}_{\alpha}=\sum_{\alpha'\in N(\alpha)} h^{(l)}_{\alpha'}.
\]
The update is coordinate-wise with learnable vectors $W^{\mathrm{self}}_{l},W^{\mathrm{neigh}}_{l}\in\mathbb{R}^{H}$:
\[
h^{(l+1)}_{\alpha}=\mathrm{ReLU}\!\left(
W^{\mathrm{self}}_{l}\odot h^{(l)}_{\alpha}
+ W^{\mathrm{neigh}}_{l}\odot m^{(l)}_{\alpha}
+ b_{l}
\right)
\]
with $b_{l}\in\mathbb{R}^{H}$ for layer $l$. Dropout is applied after each message-passing nonlinearity. Finally, the model produces a node-level probability via a logistic output layer:
\[
\mathrm{logit}_\alpha 
= \frac{1}{H}\sum_{k=1}^{H} \bigl(w_{\mathrm{out},k} \cdot h^{(L)}_{\alpha,k}\bigr) 
+ b_{\mathrm{out}}
\]
\[\hat{x}_\alpha = \sigma(\mathrm{logit}_\alpha) \in [0,1]\]
where $w_{\mathrm{out}}\in\mathbb{R}^{H}$ and $b_{\mathrm{out}}\in\mathbb{R}$. Here $\sigma(\cdot)$ denotes the sigmoid function
\[
\sigma(z) = \frac{1}{1+e^{-z}},
\]
so that $\hat{x}_\alpha$ can be interpreted as the predicted probability that key $\alpha$ belongs to the executed plan. The model is trained by minimizing the binary cross-entropy loss
\[
\mathcal{L} = -\frac{1}{K}\sum_{\alpha=1}^{K}
\Bigl[
\tilde{x}_\alpha \log \hat{x}_\alpha
+ (1-\tilde{x}_\alpha)\log(1-\hat{x}_\alpha)
\Bigr],
\]
where $\tilde{x}_\alpha \in \{0,1\}$ is the ground-truth label indicating whether key $\alpha$ is active in the executed plan, and $K$ is the total number of key-level predictions per training instance.
The architecture and training hyperparameters are fixed across all experiments, as summarized in Table~\ref{tab:gnn_config}, ensuring that any observed performance differences are attributable to the proposed learning--optimization interaction, and not to hyperparameter selection.

\begin{table}[t]
\centering
\small
\begin{tabular}{ll}
\toprule
\textbf{Component} & \textbf{Configuration} \\
\midrule
Number of message-passing layers $L$ & 2 \\
Hidden dimension $H$    & 8 \\
Activation (hidden)     & ReLU \\
Output activation       & Sigmoid \\
Input feature dimension & $R = 15$ \\
Loss function           & Binary cross-entropy \\
Optimizer               & Adam \\
Learning rate           & $10^{-3}$ \\
Epochs                  & 300 \\
Gradient clipping       & $\|\nabla\|_2 \leq 1.0$ \\
Dropout (hidden)        & $p = 0.1$ \\
\bottomrule
\end{tabular}
\caption{Graph Neural Network architecture and training configuration.}
\label{tab:gnn_config}
\end{table}
\subsubsection{Augmented optimization problem using the GNN}
\label{subsec:gnn_aop}
Once trained, the GNN forward pass is encoded as a set of linear and integer
constraints, so that the solver can jointly optimize over the original routing
and scheduling decisions and the learned feasibility predictions. Let
$\mathrm{logit}_\alpha$ denote the output of the final readout layer associated
with key $\alpha$, and let $\hat{f}^{GNN}_{\theta}(\alpha)=\sigma(\mathrm{logit}_\alpha)$
be the predicted feasibility score. For each message-passing layer $l$, neighbor aggregation is the
unweighted sum over the neighbor set $N(\alpha)$,
\[
m_\alpha^{(l)} = \sum_{\alpha' \in N(\alpha)} h_{\alpha'}^{(l)},
\]
followed by the coordinate-wise linear transformation
\[
z_\alpha^{(l+1)} = W^{\text{self}}_{l} \odot h_\alpha^{(l)}
                 + W^{\text{neigh}}_{l} \odot m_\alpha^{(l)} + b_{l},
\]
and a ReLU activation $h_\alpha^{(l+1)} = \mathrm{ReLU}(z_\alpha^{(l+1)})$ enforced
component-wise through standard big-$M$ linearization:
\[
\begin{aligned}
& h_\alpha^{(l+1)} \ge z_\alpha^{(l+1)}, \\
& h_\alpha^{(l+1)} \le z_\alpha^{(l+1)} + M_l (1 - s_\alpha^{(l+1)}), \\
& h_\alpha^{(l+1)} \le M_l\, s_\alpha^{(l+1)}, \\
& z_\alpha^{(l+1)} \ge -M_l (1 - s_\alpha^{(l+1)}), \\
& h_\alpha^{(l+1)} \ge 0,
\end{aligned}
\]
where $s_\alpha^{(l+1)} \in \{0,1\}^{H}$ are auxiliary binary variables, one per
hidden coordinate, and the inequalities hold element-wise. The readout layer
then produces the key-level logit
\[
\mathrm{logit}_\alpha
   = \frac{1}{H}\sum_{k=1}^{H} w_{\text{out},k}\, h_{\alpha,k}^{(L)} + b_{\text{out}},
\]
followed by the sigmoid activation
\[
\hat{x}_\alpha = \sigma(\mathrm{logit}_\alpha) \in [0,1].
\]
Because the sigmoid $\sigma(\cdot)$ is non-linear and therefore not directly
representable in a mixed-integer linear program, the relation
$\hat{x}_\alpha = \sigma(\mathrm{logit}_\alpha)$ is enforced in the MILP through
a \emph{piecewise-linear (PWL) approximation} of the sigmoid. The pre-activation
logit is restricted to the interval $[-S, S]$, with $S=8$ chosen so that
$\sigma(-S)\approx 0$ and $\sigma(S)\approx 1$, and the interval is partitioned
into $K$ equal segments with breakpoints
$\{\xi_0, \xi_1, \dots, \xi_K\}$ and matching ordinates
$\{\sigma(\xi_0), \sigma(\xi_1), \dots, \sigma(\xi_K)\}$. The sigmoid relation
is then imposed via Gurobi's general piecewise-linear constraint
\[
\hat{x}_\alpha \;=\; \mathrm{PWL}_{\{\xi_k,\,\sigma(\xi_k)\}_{k=0}^{K}}\!\bigl(\mathrm{logit}_\alpha\bigr).
\]
In our implementation we use $S=8$ and $K=16$ segments, which provides a maximum
pointwise approximation error below $5\times 10^{-3}$ on $[-S,S]$ while keeping
the number of auxiliary variables and constraints introduced per key small.
Together, the message-passing, linear-transformation, ReLU and PWL-sigmoid
equations form a set of linear and integer constraints that represent the full
GNN forward pass within the MILP. Substituting these constraints for the
abstract predictor identity $\hat{x}_{\alpha} = \hat{f}^{m}_{\theta}(\mathcal{G},X,{\alpha})$
in the augmented optimization problem of Section~\ref{sec:CLIF_framework}
yields the GNN-augmented model $\mathcal{LTRSP}_{\text{AUG}}^{GNN}$, while its
objective and routing and scheduling constraints remain unchanged. The solver
therefore jointly optimizes over the original decision variables and the learned
feasibility predictions within a single mixed-integer program.

\raggedbottom
\subsection{Stacking mechanism}
\label{sec:stacking}
When multiple predictors are embedded simultaneously in the AOP, a natural question arises: which prediction should the solver follow for a given decision. Since different implicit constraints may be of different natures, a single predictor may not uniformly dominate across all decisions. We address this through a confidence-based stacking mechanism that selects, for each decision, the most confident predictor among the three embedded models.

Let $\hat{x}_{\alpha}^{\mathrm{LR}}$, $\hat{x}_{\alpha}^{\mathrm{DT}}$, and $\hat{x}_{\alpha}^{\mathrm{GNN}}$ denote the feasibility predictions of the linear regressor, the decision tree, and the GNN for key $\alpha$, respectively. Each prediction lies in $[0,1]$, where values close to $0$ or $1$ indicate high confidence and values near $0.5$ indicate uncertainty. We define the confidence of each predictor as its deviation from the indecision point:
\[ \gamma_{\alpha}^{m} = \left|\hat{x}_{\alpha}^{m} - 0.5\right|, \qquad m \in \mathcal{M} = \{\mathrm{LR}, \mathrm{DT}, \mathrm{GNN}\}. \]
The stacking mechanism selects the predictor with the highest confidence. We introduce binary selection variables $s_{\alpha}^{m} \in \{0,1\}$ for each predictor $m$ and key $\alpha$, subject to
\[ s_{\alpha}^{\mathrm{LR}} + s_{\alpha}^{\mathrm{DT}} + s_{\alpha}^{\mathrm{GNN}} = 1, \qquad \forall\, \alpha \in \mathcal{A}. \]
To ensure that only the most confident predictor is selected, we introduce an auxiliary variable $\gamma_{\alpha}^{\max} \geq \gamma_{\alpha}^{m}$ for all $m$, and enforce
\[ \gamma_{\alpha}^{m} \geq \gamma_{\alpha}^{\max} - M_{\gamma}(1 - s_{\alpha}^{m}), \quad \forall m \in \mathcal{M}. \]
where $M_{\gamma}$ is a sufficiently large constant (in practice, $M_{\gamma} = 0.5$ since $\gamma_{\alpha}^{m} \in [0, 0.5]$). This ensures that $s_{\alpha}^{m} = 1$ only if predictor $m$ achieves the maximum confidence for key $\alpha$. The stacked prediction $\hat{x}_{\alpha}^{\mathrm{STK}}$ is then linked to the selected predictor through big-$M$ coupling constraints, for all $m \in \mathcal{M}$:
\[ \hat{x}_{\alpha}^{m} - M_{c}(1 - s_{\alpha}^{m}) \;\leq\; \hat{x}_{\alpha}^{\mathrm{STK}} \;\leq\; \hat{x}_{\alpha}^{m} + M_{c}(1 - s_{\alpha}^{m}). \]
where $M_{c} = 1$ since all predictions lie in $[0,1]$. When $s_{\alpha}^{m} = 1$, these constraints enforce $\hat{x}_{\alpha}^{\mathrm{STK}} = \hat{x}_{\alpha}^{m}$; when $s_{\alpha}^{m} = 0$, they become inactive. The deviation variables in the AOP are then defined with respect to the stacked prediction, for all $\alpha \in \mathcal{A}$:
\[ \Delta_{\alpha} \geq x_{\alpha} - \hat{x}_{\alpha}^{\mathrm{STK}}, \qquad \Delta_{\alpha} \geq \hat{x}_{\alpha}^{\mathrm{STK}} - x_{\alpha}. \]
The complete stacked AOP, denoted $\mathcal{LTRSP}_{\mathrm{AUG}}^{\mathrm{STK}}$, thus takes the form
\begin{longtable}{@{}r@{\;}l@{\hspace{2em}}l@{}}
$\displaystyle\min_{\mathbf{x},\,\hat{\mathbf{x}},\,\hat{\mathbf{x}}^{\mathrm{STK}},\,\boldsymbol{\Delta},\,\mathbf{s},\,\boldsymbol{\gamma}^{\max}}$
& $\displaystyle\sum_{\alpha \in \mathcal{A}} c_{\alpha} x_{\alpha} + \lambda \sum_{\alpha \in \mathcal{A}} \Delta_{\alpha}$ & \\[8pt]
s.t.
& $\mathbf{A}\,\mathbf{x} \;\geq\; \mathbf{b},\quad x_{\alpha} \in \{0,1\}$ & $\forall\, \alpha$ \\[3pt]
& $\hat{x}_{\alpha}^{m} = \hat{f}^m_\theta(\mathcal{G}, X)_{\alpha}$ & $\forall\, \alpha,\, m$ \\[3pt]
& $\gamma_{\alpha}^{m} = |\hat{x}_{\alpha}^{m} - 0.5|$ & $\forall\, \alpha,\, m$ \\[3pt]
& $\sum_{m} s_{\alpha}^{m} = 1$ & $\forall\, \alpha$ \\[3pt]
& $\gamma_{\alpha}^{\max} \geq \gamma_{\alpha}^{m}$ & $\forall\, \alpha,\, m$ \\[3pt]
& $\gamma_{\alpha}^{m} \geq \gamma_{\alpha}^{\max} - M_{\gamma}(1 - s_{\alpha}^{m})$ & $\forall\, \alpha,\, m$ \\[3pt]
& $\hat{x}_{\alpha}^{m} - M_{c}(1\!-\!s_{\alpha}^{m}) \leq \hat{x}_{\alpha}^{\mathrm{STK}}$ & $\forall\, \alpha,\, m$ \\[3pt]
& $\hat{x}_{\alpha}^{\mathrm{STK}} \leq \hat{x}_{\alpha}^{m} + M_{c}(1\!-\!s_{\alpha}^{m})$ & $\forall\, \alpha,\, m$ \\[3pt]
& $\Delta_{\alpha} \geq x_{\alpha} - \hat{x}_{\alpha}^{\mathrm{STK}}$ & $\forall\, \alpha$ \\[3pt]
& $\Delta_{\alpha} \geq \hat{x}_{\alpha}^{\mathrm{STK}} - x_{\alpha}$ & $\forall\, \alpha$ \\[3pt]
& $\Delta_{\alpha} \geq 0$ & $\forall\, \alpha$ \\[3pt]
& $s_{\alpha}^{m} \in \{0,1\}$ & $\forall\, \alpha,\, m$ \\
\end{longtable}
This formulation introduces $3|\mathcal{A}|$ additional binary variables. The solver jointly optimizes the routing and scheduling decisions and the predictor selection, effectively choosing the most reliable feasibility signal for each key within the optimization process itself. This is a key distinction from conventional ensemble methods where model selection occurs outside the optimization: here, the selection is endogenous to the MILP and benefits from the global view of the solver.
A natural question is whether embedding the three predictors jointly, rather than committing to one of them, can ever degrade the quality of the solution. The following remark shows that it cannot. For each predictor $m \in \mathcal{M}$, let $Z^{m}(\lambda)$ denote the optimal objective value of the single-predictor problem $\mathcal{LTRSP}_{\mathrm{AUG}}^{m}$ at penalty $\lambda$, and let $Z^{m}_{\mathrm{LP}}(\lambda)$ denote the optimal value of its linear relaxation, obtained by relaxing $x_{\alpha} \in \{0,1\}$ to $x_{\alpha} \in [0,1]$. Similarly, let $Z^{\mathrm{STK}}(\lambda)$ denote the optimal objective value of the stacked problem $\mathcal{LTRSP}_{\mathrm{AUG}}^{\mathrm{STK}}$, and $Z^{\mathrm{STK}}_{\mathrm{LP}}(\lambda)$ the optimal value of its linear programming relaxation, in which both $x_{\alpha}$ and $s_{\alpha}^{m}$ are relaxed to $[0,1]$.
\begin{remark}\label{rem:containment}
The stacked formulation contains each single-predictor formulation as a special case. Fixing $s_{\alpha}^{m} = 1$ for one predictor and $0$ for the others reduces the stacked deviation term to that of $\mathcal{LTRSP}_{\mathrm{AUG}}^{m}$, while the routing and scheduling constraints $\mathbf{A}\mathbf{x} \geq \mathbf{b}$ and the predictions $\hat{x}_{\alpha}^{m} = \hat{f}^m_\theta(\mathcal{G}, X)_{\alpha}$ are shared across all formulations. Hence, for every $\lambda \geq 0$,
\[ Z^{\mathrm{STK}}_{\mathrm{LP}}(\lambda)\!\leq\!\min\nolimits_{m} Z^{m}_{\mathrm{LP}}(\lambda),\; Z^{\mathrm{STK}}(\lambda)\!\leq\!\min\nolimits_{m} Z^{m}(\lambda), \]
both at the linear relaxation and at the integer level. A formal proof is given in Appendix~\ref{app:containment_proof}.
\end{remark}
By embedding the three predictors jointly and letting the solver select among them endogenously, the stacked formulation cannot underperform the best individual predictor, while retaining the freedom to combine them. The empirical magnitude of this improvement is reported in Section~\ref{sec:experiments}.
\section{Experimentation}
\label{sec:experiments}
In this section, we evaluate the proposed framework \textbf{CLIF} on a large-scale and real-world routing dataset of $\mathcal{LTRSP}$ obtained from an 
industrial forestry partner in Canada. The objective is to assess 
whether the learned predictors improve implicit operational 
feasibility, generalization to new contexts, and overall optimization 
quality, both when used as single embedded predictors and when combined through the proposed stacking mechanism. In particular, since our framework builds upon the methodological lineage of \cite{hewitt2020data} and \cite{bayani2024learning}, we explicitly include their approaches as baselines in our experimental comparison: the linear regression embedding $\mathcal{LTRSP}_{\text{AUG}}^{\text{LR}}$ corresponds to the data-driven model customization framework of \cite{hewitt2020data}, while the decision tree embedding $\mathcal{LTRSP}_{\text{AUG}}^{\text{DT}}$ corresponds to the implicit constraint learning approach of \cite{bayani2024learning}. The GNN embedding $\mathcal{LTRSP}_{\text{AUG}}^{\text{GNN}}$ and the stacked model $\mathcal{LTRSP}_{\text{AUG}}^{\text{STK}}$ are the two contributions of the present work, and are evaluated against these two baselines from the literature. We introduce first characteristics of the routing dataset and we shed light on the implicit constraints generation process. Then, we present numerical results, and we finalize this section by a discussion and insights.
All computational experiments were conducted on a workstation running
AlmaLinux~9.7, equipped with two Intel\textsuperscript{\textregistered}
Xeon\textsuperscript{\textregistered} E5-2650~v4 processors operating at 2.20~GHz (12 cores per socket, 24 cores total), 64~GB of system memory, and using a Linux~5.14.0 kernel. The Integer Linear Programs (ILPs) were solved using Gurobi Optimization (version~12.0.0). The GNN
models are implemented using the PyTorch library.
\subsection{Data description}
The dataset consists of $48$ months of weekly routing and scheduling decisions collected from a large-scale industrial forest partner. Each instance corresponds to a full operational week and is defined over a routing and scheduling graph whose nodes are selected from a pool of more than $900$ forest blocks and $50$ mills, and whose arcs encode feasible transportation links for a heterogeneous fleet of approximately $110$ trucks. From this operational history, we extract a set $\mathcal{I}$ of $192$ representative weekly instances.

For each instance $i \in \mathcal{I}$, we store two solutions defined over the same routing graph: an optimal solution $\mathbf{x}_i^\star$ of the corresponding $\mathcal{LTRSP}$, obtained by solving the formulation of Appendix~\ref{LTRSP}, and an executed routing and scheduling plan $\hat{\mathbf{x}}_i$.

In practice, for our forest partner, planners do not always log their modifications: they adjust the plan continuously, without separating recurrent operational rules from reactions to disruptions such as demand fluctuations, road closures, or equipment failures. Observed executed plans therefore mix the deterministic rules we seek to learn with these stochastic, non-recurrent disruptions, and cannot serve directly as a learning signal. We therefore construct surrogate adjusted plans that imitate planner behavior under implicit feasibility constraints: each rule is enforced on the optimal plan through local branching, while stochastic disruptions remain out of scope.

Formally, $\hat{\mathbf{x}}_i$ is obtained by solving a local branching problem centered at $\mathbf{x}_i^\star$, in which a predefined subset of rules is imposed as additional constraints; this keeps the executed plan close to the optimal one while enforcing the selected rules. The local branching construction is detailed in Appendix~\ref{localbranching}, and the full set of operational rules is reported in Table~\ref{tab:rules}.

The operational rules $R_1, \ldots, R_5$ were fixed following discussions with the planners at the industrial partner, and correspond to the categories of implicit constraints they report applying most frequently in daily practice. For each rule $R_k$, we add to $\mathcal{LTRSP}_G$ the linear constraints that encode $R_k$ and re-solve the model, so that the resulting plan satisfies $R_k$ while remaining feasible for the original routing and scheduling constraints. The executed plan $\hat{\mathbf{x}}_i$ is therefore the cost-minimal solution consistent with $R_k$ for instance $i$. The resulting plans were reviewed by the planners of the industrial partner, who confirmed that they constitute realistic and operationally valid executed plans for the corresponding rules.

The resulting paired samples $\{(\mathbf{x}_i^\star, \hat{\mathbf{x}}_i)\}_{i \in \mathcal{I}}$ are partitioned chronologically into three disjoint subsets: a training set $\mathcal{I}^{\text{train}}$ of $150$ weekly instances, a validation set $\mathcal{I}^{\text{validation}}$ of $20$ instances used for hyperparameter selection, and a test set $\mathcal{I}^{\text{test}}$ of $20$ instances reserved for the final evaluation and described in Appendix~\ref{appendix:datasets}. We denote the training, validation, and test datasets as follows :
\[
\mathcal{D}^{s} = \{\mathcal{D}_i\}_{i \in \mathcal{I}^{s}}, \quad s \in \{\text{train}, \text{validation}, \text{test}\}.
\]
The chronological split ensures that no information from the validation or test periods leaks into training, mimicking the deployment setting where the system must generalize to future operational weeks.

\begin{table*}
\centering
\caption{Operational rules settings for executed-plan generation. Each decision variable $x_{\alpha}$ corresponds to a key $\alpha = (a, v, t)$, where $a = (u, w)$ is an arc from origin node $u$ to destination node $w$, $v$ is a vehicle, and $t$ is a time index. Each node is written $n = (\ell, \tau)$ where $\ell$ denotes the physical location (homebase, block, or mill) and $\tau$ the corresponding time interval.}
\label{tab:rules}
\small
\renewcommand{\arraystretch}{1.6}
\setlength{\tabcolsep}{6pt}
\begin{tabularx}{\textwidth}{>{\raggedright\arraybackslash}p{5.5cm} >{\raggedright\arraybackslash}X}
\toprule
\textbf{Rule} & \textbf{Mathematical representation} \\
\midrule
$\mathbf{R_1}$: Vehicles $V^* \subseteq V$ cannot visit block $f^*$
&
$x_{\alpha} = 0, \quad \forall\, \alpha = (a, v, t) \text{ with } a = (u, w),\ u = (\ell_u, \tau_u),\ w = (\ell_w, \tau_w) \text{ s.t. } f^* \in \{\ell_u, \ell_w\},\ v \in V^*$
\\
\midrule
$\mathbf{R_2}$: Homebase $h^*$ cannot dispatch to blocks $F^* \subseteq F$
&
$x_{\alpha} = 0, \quad \forall\, \alpha = (a, v, t) \text{ with } a = (u, w),\ u = (h^*, \tau_u),\ w = (f, \tau_w) \text{ s.t. } f \in F^*$
\\
\midrule
$\mathbf{R_3}$: Vehicle group $V^*$ cannot visit block groups $F_A$ and $F_B$ on the same day
&
$\displaystyle
\zeta_A = \mathbf{1}\!\Bigg[\sum_{\substack{\alpha = (a,v,t):\\ a = (u,w),\, u = (\ell_u,\tau_u),\, w = (\ell_w,\tau_w),\\ \{\ell_u,\ell_w\} \cap F_A \neq \emptyset,\, v \in V^*}} x_{\alpha} \ge 1\Bigg]$
\\
& $\displaystyle
\zeta_B = \mathbf{1}\!\Bigg[\sum_{\substack{\alpha = (a,v,t):\\ a = (u,w),\, u = (\ell_u,\tau_u),\, w = (\ell_w,\tau_w),\\ \{\ell_u,\ell_w\} \cap F_B \neq \emptyset,\, v \in V^*}} x_{\alpha} \ge 1\Bigg]$
\\
& $\zeta_A + \zeta_B \le 1$
\\
\midrule
$\mathbf{R_4}$: Mills $M^* \subseteq M$ cannot be served after time interval $j^*$
&
$x_{\alpha} = 0, \quad \forall\, \alpha = (a, v, t) \text{ with } a = (u, w),\ w = (m, j) \text{ s.t. } m \in M^*,\ j > j^*$
\\
\midrule
$\mathbf{R_5}$: At most one of the mill groups $M_1, M_2 \subseteq M$ can be served in the morning ($j < j^*$)
&
$\displaystyle
\zeta_k = \mathbf{1}\!\Bigg[\sum_{\substack{\alpha = (a,v,t):\\ a = (u,w),\, w = (m, j),\\ m \in M_k,\, j < j^*}} x_{\alpha} \ge 1\Bigg],\ k \in \{1, 2\},\ \zeta_1 + \zeta_2 \le 1$
\\
\bottomrule
\end{tabularx}
\end{table*}

\subsection{Numerical results}
\begin{figure}[htbp]
    \centering
    \includegraphics[width=0.70\textwidth]{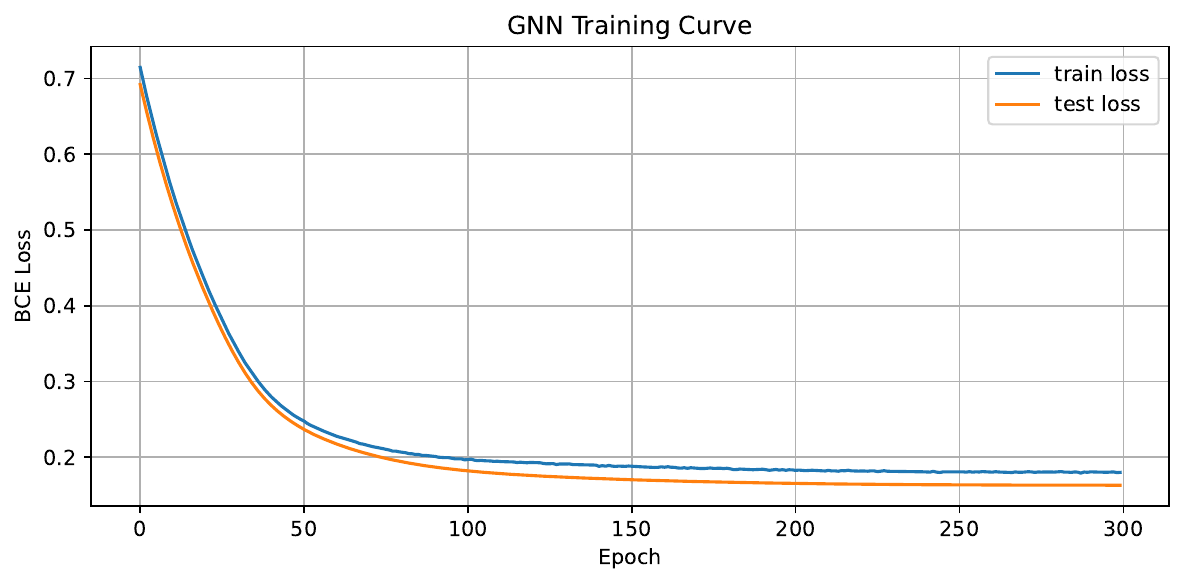}
    \caption{\centering GNN binary cross entropy curve}
    \label{fig:gnn curve}
\end{figure}

We use the training dataset $\mathcal{D}^{\text{train}}$ to learn the parameters
of all three predictors: the linear regressor $\hat{f}_{\theta}^{LR}$ following the formulation of \cite{hewitt2020data}, the decision tree $\hat{f}_{\theta}^{DT}$ following the formulation of \cite{bayani2024learning}, and the GNN $\hat{f}_{\theta}^{GNN}$ proposed in this work. For the GNN, we additionally use $\mathcal{D}^{\text{validation}}$ as a validation set to monitor convergence during training.
Figure~\ref{fig:gnn curve} reports the evolution of the binary cross-entropy loss on both the training and validation sets throughout the learning process when the executed plans mimic rule $R_1$. The curves exhibit a smooth and monotonic decrease, indicating stable convergence of the learning procedure. Moreover, the absence of a significant gap between training and validation losses suggests that the model generalizes well and does not suffer from overfitting. The loss stabilizes after approximately 200 epochs, confirming that the learned node representations have reached a consistent predictive regime.
The performance
of the proposed framework is then evaluated by solving $\mathcal{LTRSP}_{\mathrm{AUG}}$ for each test instance
$i \in \mathcal{I}^{\text{test}}$. 

To evaluate the quality of the obtained solutions, we rely on two performance metrics, namely the satisfaction rate and the objective gap. These
metrics respectively measure the ability of $\mathcal{LTRSP}_{\mathrm{AUG}}$ to recover feasible
solutions with respect to operational implicit  rules, and the degradation in objective value induced by enforcing such rules through soft constraints learned by \textbf{CLIF}.

The satisfaction rate measures the proportion of test instances for which the
solution returned by the $\mathcal{LTRSP}_{\mathrm{AUG}}$  satisfies all hidden operational rules. Let
$\mathbf{x}_i^{\text{AUG}}$ denote the solution obtained by solving the $\mathcal{LTRSP}_{\mathrm{AUG}}$ for instance
$i$, and let $R \in \mathcal{R} = \{R_1, R_2, \ldots, R_5\}$ denote the set of applied
business rules. We define the indicator variable
\[
\mathbb{I}_i =
\begin{cases}
    1, & \text{if } \mathbf{x}_i^{\text{AUG}} \text{ satisfies activated rules in } R, \\
0, & \text{otherwise}.
\end{cases}
\]
The satisfaction rate is then computed as
\[
\textit{SatRate} =
\frac{1}{|\mathcal{I}^{\text{test}}|}
\sum_{i \in \mathcal{I}^{\text{test}}} \mathbb{I}_i.
\]

While enforcing hidden rules enhances feasibility, it can also deteriorate the objective value. To quantify this trade-off, we measure the relative objective gap between $Z_i^{\mathcal{LTRSP}_G\text{-}\mathcal{IC}}$, the objective value of $\mathcal{LTRSP}$ for instance $i$ when solved with the complete set of implicit constraints assumed to be fully known in advance, and $Z_i^{\mathcal{LTRSP}_{\mathrm{AUG}}}$, the objective value obtained by $\mathcal{LTRSP}_{\mathrm{AUG}}$ for the same instance.
The relative gap for instance $i$ is defined as
\begin{equation}
\textit{Gap}_i =
\frac{
Z_i^{\mathcal{LTRSP}_G\text{-}\mathcal{IC}} - Z_i^{\mathcal{LTRSP}_{\mathrm{AUG}}}
}{
Z_i^{\mathcal{LTRSP}_G\text{-}\mathcal{IC}}
}.
\end{equation}
Since this measure is meaningful only for instances that satisfy all concealed
rules, we report the average gap over all feasible test instances:
\[
\overline{\textit{Gap}} =
\frac{1}{\sum_{i \in \mathcal{I}^{\text{test}}} \mathbb{I}_i}
\sum_{i \in \mathcal{I}^{\text{test}}} \mathbb{I}_i \cdot \textit{Gap}_i.
\]

\begin{figure*}[!tbp]
    \centering
    \begin{subfigure}[t]{0.49\textwidth}
        \centering
        \includegraphics[width=\linewidth]{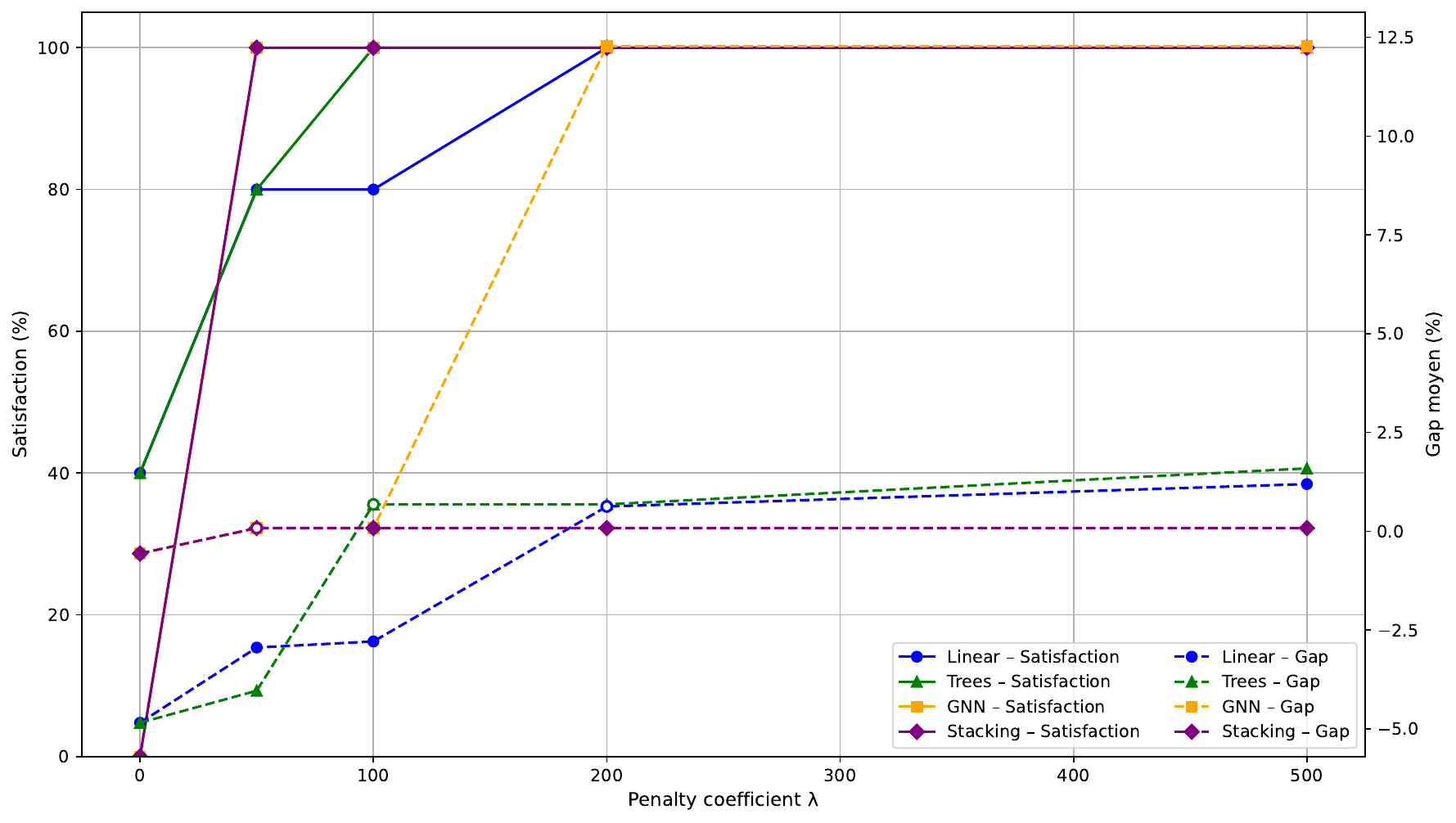}
        \caption{Rule $R_1$}
        \label{fig:rule 1}
    \end{subfigure}
    \hfill
    \begin{subfigure}[t]{0.49\textwidth}
        \centering
        \includegraphics[width=\linewidth]{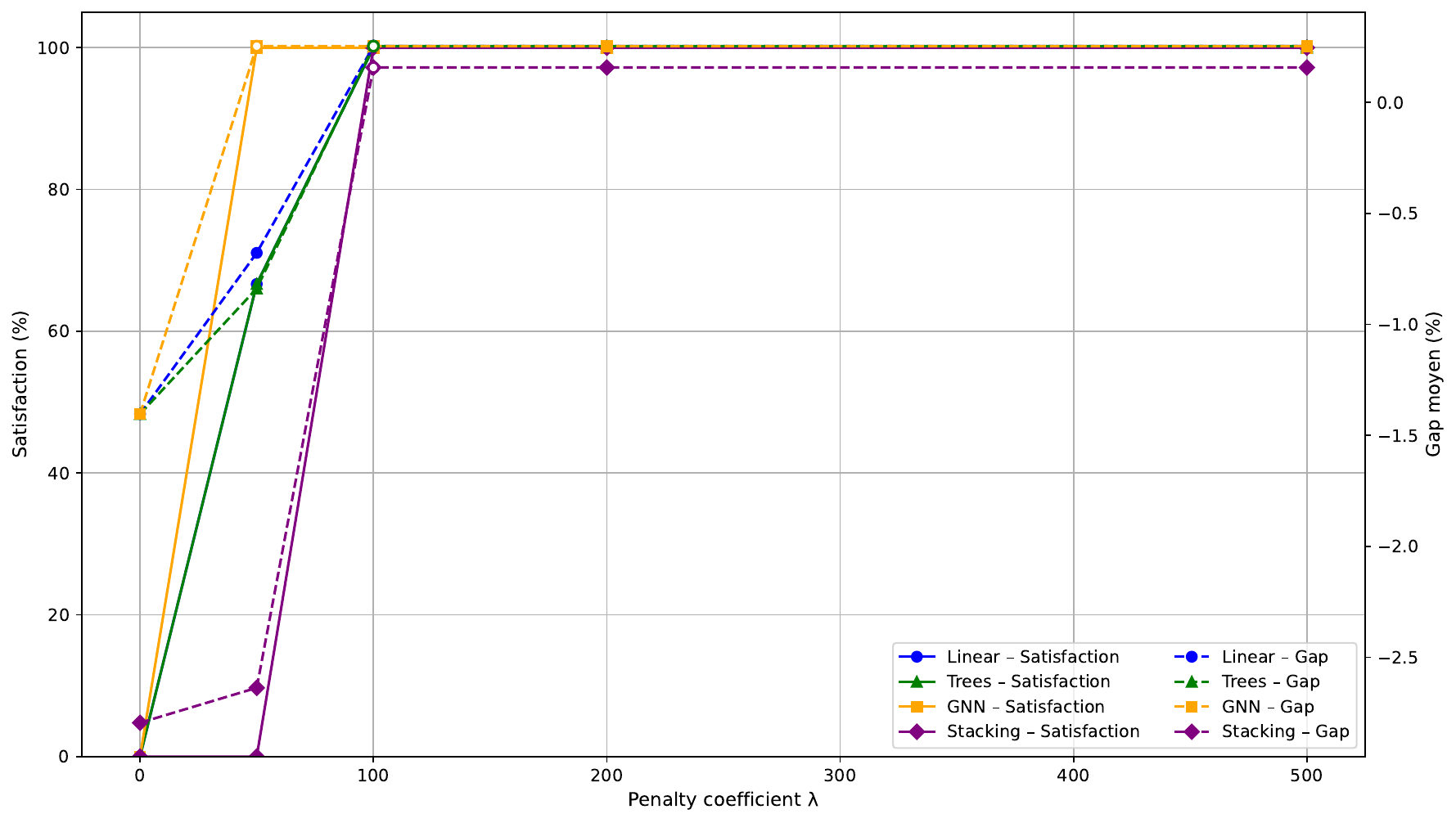}
        \caption{Rule $R_2$}
        \label{fig:rule 2}
    \end{subfigure}

    \vspace{0.8em}

    \begin{subfigure}[t]{0.49\textwidth}
        \centering
        \includegraphics[width=\linewidth]{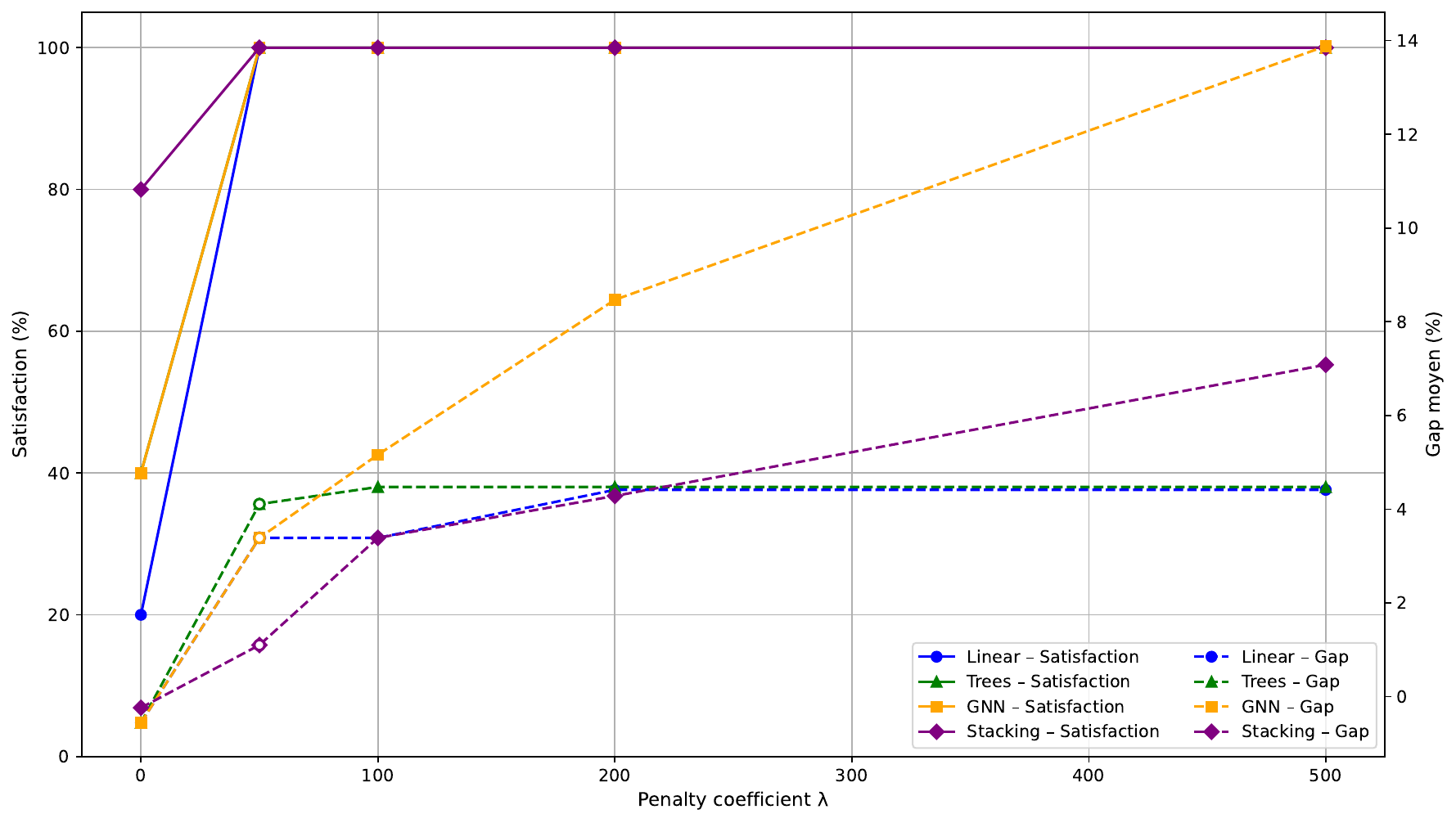}
        \caption{Rule $R_3$}
        \label{fig:rule 3}
    \end{subfigure}
    \hfill
    \begin{subfigure}[t]{0.49\textwidth}
        \centering
        \includegraphics[width=\linewidth]{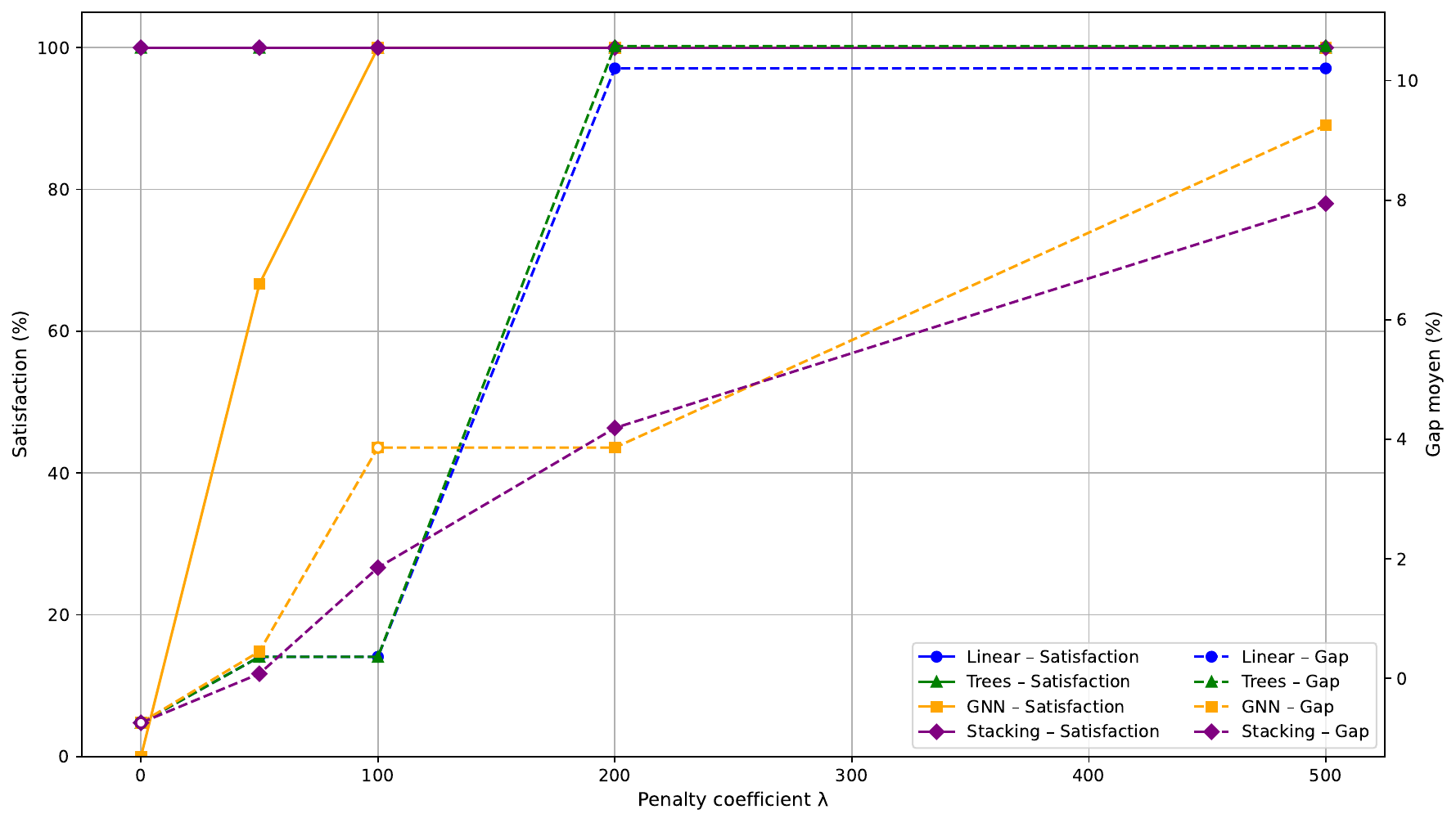}
        \caption{Rule $R_4$}
        \label{fig:rule 4}
    \end{subfigure}

    \vspace{0.8em}

    \begin{subfigure}[t]{0.49\textwidth}
        \centering
        \includegraphics[width=\linewidth]{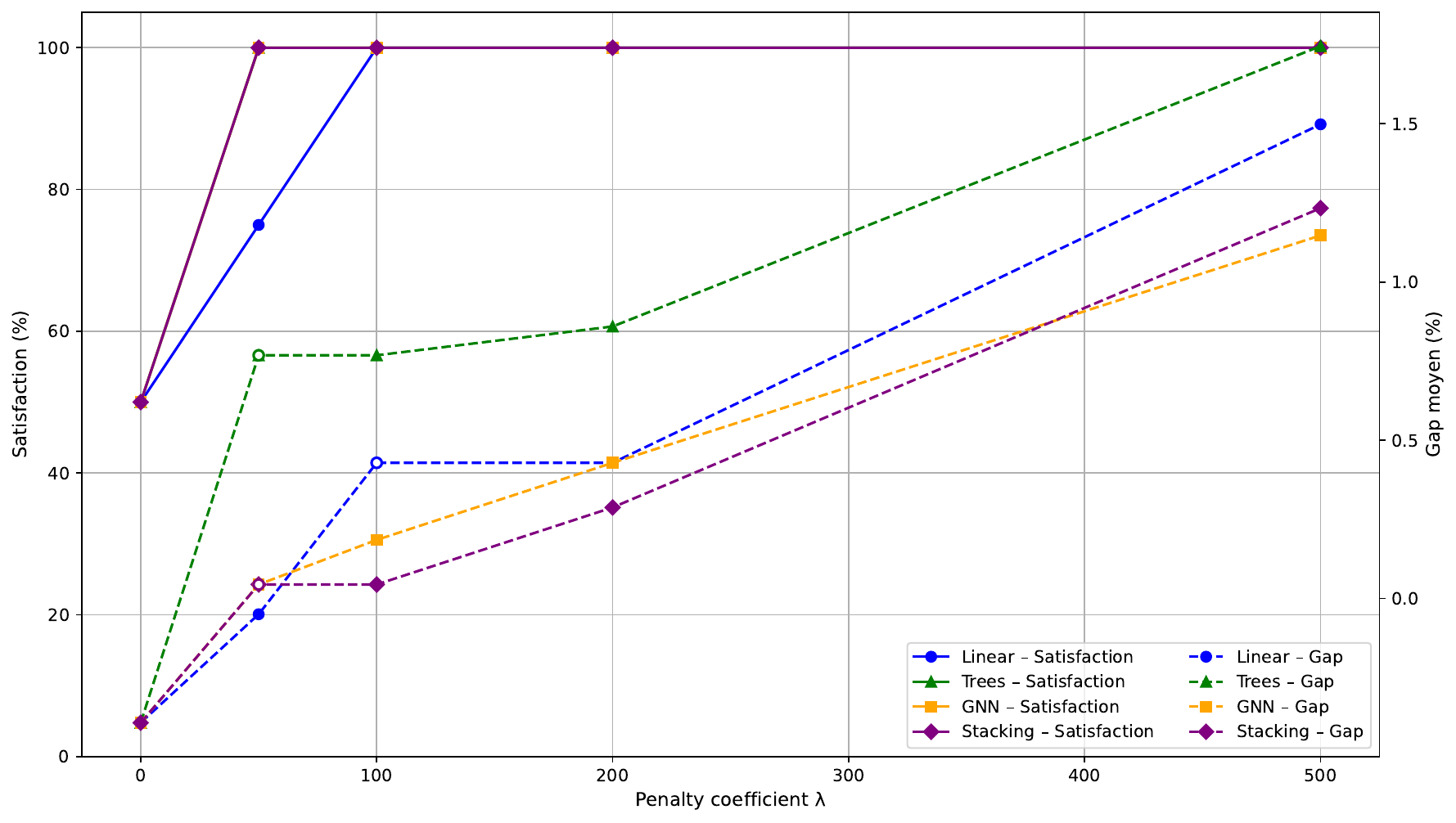}
        \caption{Rule $R_5$}
        \label{fig:rule 5}
    \end{subfigure}
    \hfill
    \begin{subfigure}[t]{0.49\textwidth}
        \centering
        \includegraphics[width=\linewidth]{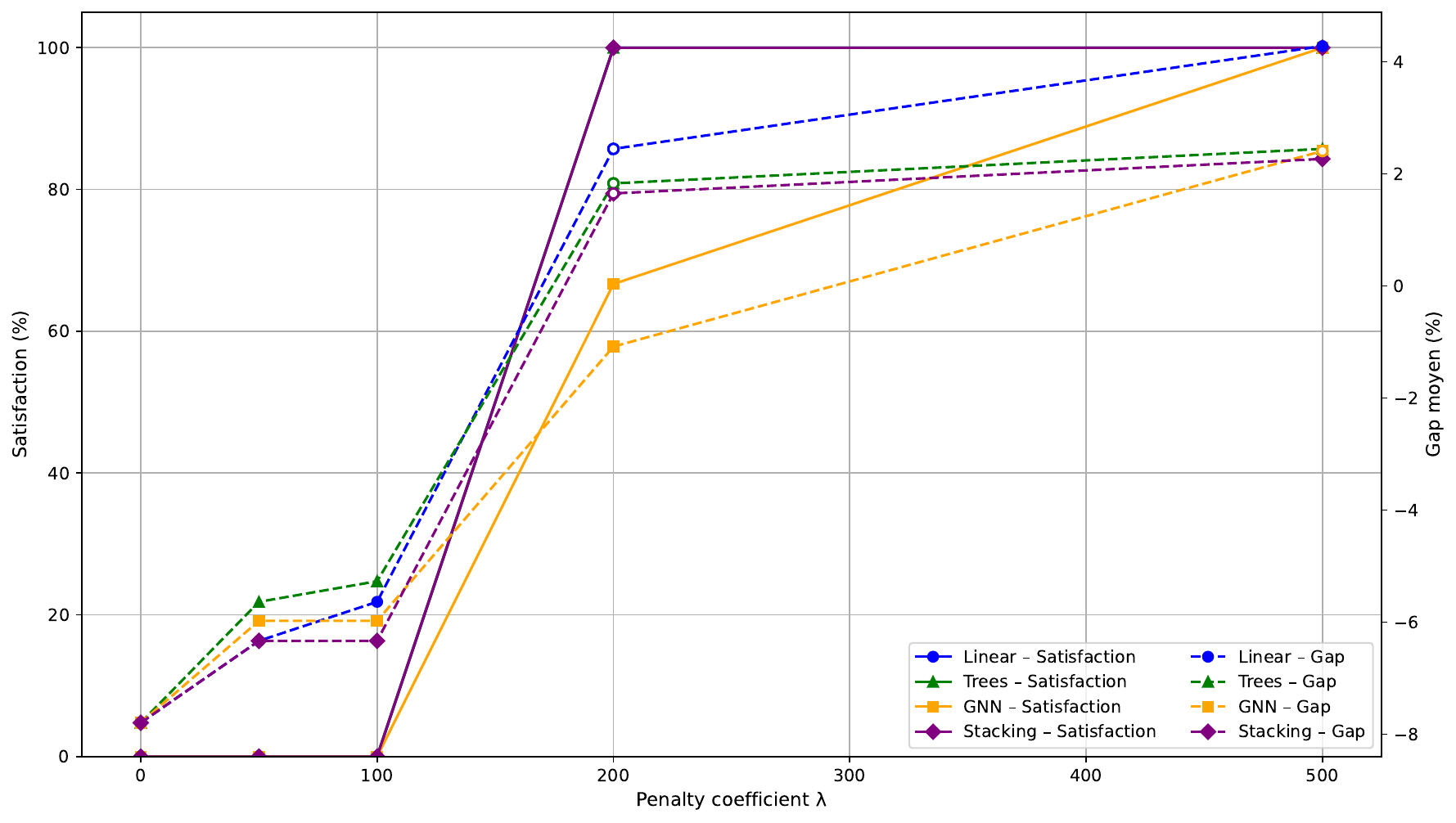}
        \caption{Combined $\{R_1,\dots,R_5\}$}
        \label{fig:rule 6}
    \end{subfigure}

    \caption{Comparison of \textbf{CLIF} variants ($\mathcal{LTRSP}_{\text{AUG}}^{\text{GNN}}$ and $\mathcal{LTRSP}_{\text{AUG}}^{\text{STK}}$) against the literature baselines of \cite{hewitt2020data} ($\mathcal{LTRSP}_{\text{AUG}}^{\text{LR}}$) and \cite{bayani2024learning} ($\mathcal{LTRSP}_{\text{AUG}}^{\text{DT}}$). Each subfigure reports the satisfaction rate (left axis, solid lines) and the objective gap (right axis, dashed lines) as functions of the penalty coefficient $\lambda$, for the five operational rules and the combined regime $\{R_1,\dots,R_5\}$.}
    \label{fig:all_rules}
\end{figure*}

We conduct six experiments, each corresponding either to an individual operational rule 
$R_k \in \{R_1, \ldots, R_5\}$ or to the combination of all five rules simultaneously. In each experiment, we vary the penalty 
coefficient $\lambda$ over a predefined range and evaluate four variants of 
$\mathcal{LTRSP}_{\text{AUG}}$ that fall into two categories. The first category gathers the two baselines from the literature: the linear regression embedding $\mathcal{LTRSP}_{\text{AUG}}^{\text{LR}}$, which reproduces the data-driven optimization model customization approach of \cite{hewitt2020data}, and the decision tree embedding $\mathcal{LTRSP}_{\text{AUG}}^{\text{DT}}$, which reproduces the implicit constraint learning approach of \cite{bayani2024learning}. The second category gathers the two contributions of the present work: the graph neural network embedding $\mathcal{LTRSP}_{\text{AUG}}^{\text{GNN}}$, which exploits the relational structure of the routing graph through message passing, and the stacked variant $\mathcal{LTRSP}_{\text{AUG}}^{\text{STK}}$, which combines all three predictors through the stacking mechanism described in Section~\ref{sec:stacking}. Numerical results are reported in Table~\ref{tab:results_comparison}, and the 
corresponding curves in Figures~\ref{fig:rule 1}--\ref{fig:rule 6} illustrate 
how both the satisfaction rate and the objective gap evolve as functions of 
$\lambda$.

Across all six experiments and all four model variants, 
$\mathcal{LTRSP}_{\text{AUG}}$ successfully propagates the learned 
feasibility signal into the solver: every variant reaches a satisfaction 
rate of $100\%$ on every rule. This uniformity follows by construction: once $\lambda$ is large enough relative
to the routing cost scale, any predictor better than random makes violations
strictly suboptimal, so the satisfaction rate saturates at $100\%$. Full
compliance is therefore guaranteed by the penalty mechanism rather than
indicative of model quality. The discriminating metric is the objective gap,
which measures the additional routing cost incurred to maintain compliance and
thus reflects how faithfully each predictor approximates the true feasible
region.

As reported in Table~\ref{tab:results_comparison}, the GNN embedding $\mathcal{LTRSP}_{\text{AUG}}^{\text{GNN}}$ matches or outperforms both baselines on every rule, illustrating the value of exploiting the routing graph topology through message passing rather than relying on feature-based predictors alone. More importantly, the stacked variant $\mathcal{LTRSP}_{\text{AUG}}^{\text{STK}}$ consistently achieves the lowest objective gap across all rules, including those for which one of the baselines was already competitive. On rules where one individual predictor already dominates, the stacked variant preserves that performance, confirming that the aggregation does not degrade the best available signal. On rules where the individual models exhibit more heterogeneous behavior, the stacked variant achieves a gap strictly lower than any single model, including the baselines from \cite{hewitt2020data} and \cite{bayani2024learning}, indicating that it effectively combines complementary information rather than simply averaging predictions.
\begin{table*}[h!]
\centering
\caption{Performance comparison of \textbf{CLIF} against the two baselines from the literature: the linear regression embedding of \cite{hewitt2020data} and the decision tree embedding of \cite{bayani2024learning}. All values are averaged over the test set under different rule settings.}
\label{tab:results_comparison}
\small
\renewcommand{\arraystretch}{1.2}
\setlength{\tabcolsep}{6pt}
\begin{tabular}{lcc|cc|cc|cc}
\toprule
& \multicolumn{4}{c|}{\textbf{Baselines from the literature}} & \multicolumn{4}{c}{\textbf{CLIF }} \\
\cmidrule(lr){2-5} \cmidrule(lr){6-9}
\multirow{2}{*}{Rule}&
\multicolumn{2}{c|}{$\mathcal{LTRSP}_{\text{AUG}}^{\text{LR}}$} &
\multicolumn{2}{c|}{$\mathcal{LTRSP}_{\text{AUG}}^{\text{DT}}$} &
\multicolumn{2}{c|}{$\mathcal{LTRSP}_{\text{AUG}}^{\text{GNN}}$} &
\multicolumn{2}{c}{$\mathcal{LTRSP}_{\text{AUG}}^{\text{STK}}$} \\
& \multicolumn{2}{c|}{\citep{hewitt2020data}} &
\multicolumn{2}{c|}{\citep{bayani2024learning}} &
\multicolumn{2}{c|}{} & \multicolumn{2}{c}{} \\
\cmidrule(lr){2-3}
\cmidrule(lr){4-5}
\cmidrule(lr){6-7}
\cmidrule(lr){8-9}
& Sat.\ (\%) & Gap (\%)
& Sat.\ (\%) & Gap (\%)
& Sat.\ (\%) & Gap (\%)
& Sat.\ (\%) & Gap (\%) \\
\midrule
$R_1$ & $100.0$ & $0.63$ & $100.0$ & $0.68$ & $100.0$ & $0.08$ & $100.0$ & $\mathbf{0.08}$ \\
$R_2$ & $100.0$ & $0.61$ & $100.0$ & $0.54$ & $100.0$ & $0.54$ & $100.0$ & $\mathbf{0.15}$ \\
$R_3$ & $100.0$ & $3.39$ & $100.0$ & $4.11$ & $100.0$ & $3.39$ & $100.0$ & $\mathbf{1.33}$ \\
$R_4$ & $100.0$ & $0.41$ & $100.0$ & $0.41$ & $100.0$ & $0.08$ & $100.0$ & $\mathbf{0.07}$ \\
$R_5$ & $100.0$ & $0.43$ & $100.0$ & $0.77$ & $100.0$ & $0.04$ & $100.0$ & $\mathbf{0.04}$ \\
\midrule
$\{R_1,\dots,R_5\}$ & $100.0$ & $2.44$ & $100.0$ & $1.82$ & $100.0$ & $2.34$ & $100.0$ & $\mathbf{1.64}$ \\
\bottomrule
\end{tabular}
\end{table*}

The combined experiment $\{R_1, \ldots, R_5\}$, in which all five rules 
must be satisfied simultaneously, further amplifies this effect. 
Enforcing multiple implicit constraints at once increases the tension 
between compliance and solution quality, as the solver must navigate a 
more restricted feasible region. In this challenging setting, both baselines exhibit a substantial degradation in objective quality, while the stacked variant maintains the lowest gap, showing that the complementarity between model families becomes 
increasingly valuable as the complexity of the implicit constraint 
landscape grows.

These results reveal a structural insight: no single predictor family 
uniformly dominates across all constraint types. The GNN captures 
relational dependencies through message passing, the decision tree 
of \cite{bayani2024learning} encodes conditional logic, and the linear regressor of \cite{hewitt2020data} estimates global 
trends, yet each exhibits weaknesses on constraint structures that 
do not align with its inductive bias. The stacking mechanism mitigates 
this limitation by allowing the solver to leverage the most faithful 
predictor for each portion of the problem, producing a closer 
approximation of the true feasible region with minimal distortion of 
the original objective.
\begin{figure}[!htbp]
  \centering
  \includegraphics[width=\linewidth]{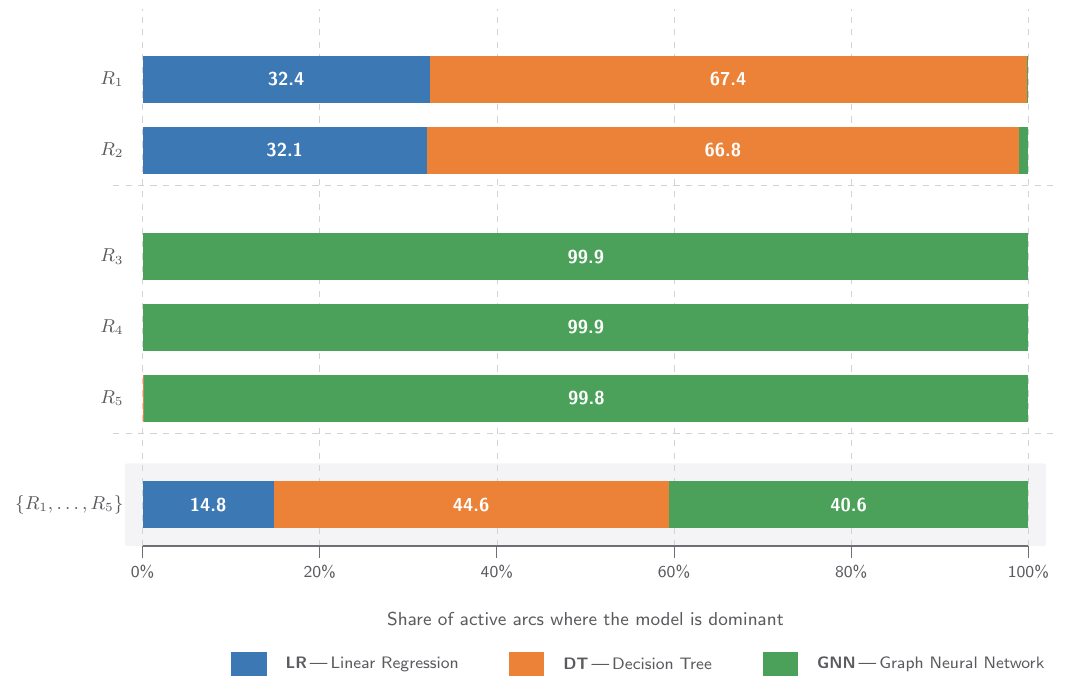}
  \caption{\centering Per-rule model dominance in the stacking approach.}
  \label{fig:dominance}
\end{figure}

Figure~\ref{fig:dominance} makes this structural insight visually 
explicit. For each rule, we report the share of active arcs on which 
each base predictor is selected as the dominant signal by the 
confidence-based stacking model. Rules $R_1$ and $R_2$, which encode 
local arc-blocking restrictions, are almost exclusively handled by the 
decision tree and linear regression, whereas the relational rules $R_3$, $R_4$ and $R_5$ are 
captured by the GNN on more than $99\%$ of active 
arcs. In the combined experiment, no individual predictor reaches more 
than $45\%$ dominance, and the three model families coexist with 
non-trivial shares. 

From a practical standpoint, this finding matters because the nature 
of implicit rules is rarely known in advance. Committing to a single 
model family, as in \cite{hewitt2020data} and \cite{bayani2024learning}, 
risks a mismatch between the predictor's modelling assumptions and 
the actual structure of the constraints. By combining predictors of 
different natures, the stacking approach proposed here removes that 
commitment and remains robust when implicit rules evolve over time.

Beyond the confidence-based aggregation retained in $\mathcal{LTRSP}_{\text{AUG}}^{\text{STK}}$, we also investigated a majority-voting layer added on top of it. In this variant, the confidence-based mechanism still selects the most confident predictor for each arc, but each of the three predictors additionally casts a binary vote, set to one if its prediction is above $0.5$ and zero otherwise. Whenever the vote of the confidence-selected model disagrees with the majority vote of the three, the majority overrides the confidence-based output; otherwise, the confidence-based output is kept unchanged. The full MILP formulation of this layered mechanism is provided in Appendix~\ref{app:voting}.

The voting layer produces a mixed empirical pattern. It improves the gap on the combined experiment $\{R_1, \ldots, R_5\}$, while degrading it on rules $R_3$ and $R_5$ and leaving the others essentially unchanged. This contrast is interpretable: when a single predictor dominates and the other two produce nearly uninformative outputs near $0.5$, turning those outputs into binary votes amplifies noise rather than reducing it. Because operational deployment requires consistent per-rule performance rather than gains conditional on the structure of the implicit constraint, we retain confidence-based stacking as the reference variant in $\mathcal{LTRSP}_{\text{AUG}}^{\text{STK}}$. Adding a majority-voting layer remains a promising direction for settings involving many simultaneous rules and could be combined with confidence-based aggregation through weighted or threshold-based mechanisms, which we leave to future work.

Beyond the mathematical performance, the augmented plans generated under each rule were submitted to the operational planners of the industrial partner for review. The planners confirmed the operational validity of the proposed plans, indicating that the rule-compliance achieved by $\mathcal{LTRSP}_{\mathrm{AUG}}$ corresponds to plans that are feasible and acceptable for execution in practice. This practitioner-level validation complements the quantitative metrics of Table~\ref{tab:results_comparison} and supports the interpretation of $\mathcal{LTRSP}_{\mathrm{AUG}}^{\text{STK}}$ as an operationally deployable framework.

\subsection{Computational analysis}
\label{sec:computational_analysis}

The practical value of \textbf{CLIF} depends not only on the predictive accuracy of the embedded models but also on the tractability of the resulting augmented optimization problem, since the system is intended to support daily operational planning. Embedding a learned predictor inside the MILP introduces additional variables and constraints, and modifies the structure of the linear programming relaxation. The relevant question is whether this additional cost remains compatible with the operational use of the system. Table~\ref{tab:computational_footprint} summarizes the results that characterize this trade-off.

\begin{table*}[h!]
\centering
\caption{Computational results of each $\mathcal{LTRSP}_{\text{AUG}}$ variant compared to the baseline $\mathcal{LTRSP}_G$. Counts are additional variables (binary and continuous) and constraints introduced by the embedding, in closed form in $|\mathcal{A}|$, $H$, and $n_{\text{leaf}}$, where $|\mathcal{A}|$ is the number of keys, $H$ the GNN hidden dimension, and $n_{\text{leaf}}$ the number of decision-tree leaves. Times are averaged over the test set under $\{R_1,\dots,R_5\}$ ($\pm$ standard deviation). The LP gap is the relative difference between the integer optimum and the linear programming relaxation value.}
\label{tab:computational_footprint}
\footnotesize
\renewcommand{\arraystretch}{1.4}
\setlength{\tabcolsep}{8pt}
\resizebox{\textwidth}{!}{%
\begin{tabular}{lccccc}
\toprule
\textbf{Variant} & \textbf{Bin.\ vars.} & \textbf{Cont.\ vars.} & \textbf{Constraints} & \textbf{Avg.\ time (min)} & \textbf{LP gap (\%)} \\
\midrule
$\mathcal{LTRSP}_G$ (baseline)              & --               & --          & --                & $10.2\;(\pm 2.1)$   & $6.08\;(\pm 5.99)$ \\
$\mathcal{LTRSP}_{\text{AUG}}^{\text{LR}}$  & $0$              & $2|\mathcal{A}|$      & $3|\mathcal{A}|$            & $10.8\;(\pm 2.4)$   & $1.88\;(\pm 1.28)$   \\
$\mathcal{LTRSP}_{\text{AUG}}^{\text{DT}}$  & $|\mathcal{A}|\cdot n_{\text{leaf}}$    & $2|\mathcal{A}|$      & $|\mathcal{A}|(4 n_{\text{leaf}}+3)$     & $12.6\;(\pm 2.9)$   & $1.93\;(\pm 1.35)$   \\
$\mathcal{LTRSP}_{\text{AUG}}^{\text{GNN}}$ & $|\mathcal{A}|\cdot H$     & $|\mathcal{A}|(2H+1)$ & $|\mathcal{A}|(5H+3)$       & $15.1\;(\pm 3.7)$  & $2.25\;(\pm 1.22)$   \\
$\mathcal{LTRSP}_{\text{AUG}}^{\text{STK}}$ & $|\mathcal{A}|(H+n_{\text{leaf}}+3)$   & $|\mathcal{A}|(2H+5)$ & $|\mathcal{A}|(5H+4 n_{\text{leaf}}+11)$ & $17.4\;(\pm 4.1)$  & $2.00\;(\pm 1.43)$   \\
\bottomrule
\end{tabular}%
}
\end{table*}

Larger models do not entail proportionally longer solving times. The linear-regression embedding is virtually free, adding only about $0.6$~minutes over the baseline despite introducing $2|\mathcal{A}|$ continuous variables and $3|\mathcal{A}|$ deviation constraints. The graph neural network embedding is substantially larger than the decision-tree embedding, yet its solving time exceeds that of the decision tree by only about $2.5$~minutes, far less than what the difference in structural size would suggest. The big-$M$ coefficients used to linearize the ReLU activations are derived from the magnitudes of the trained weights, which keeps the linearization tight and the LP relaxation informative. The leaf-selection binaries of the decision tree are coupled through a disjunctive constraint imposing that exactly one leaf is active, whose continuous relaxation is structurally looser; the discrete structure of the tree is therefore revealed mainly through branching. The combinatorial cost of embedding a neural predictor is thus not proportional to the count of binaries it introduces. The same logic extends to the stacked formulation: although it is the largest model on every structural dimension, its solving time remains close to that of the GNN alone, exceeding it by roughly $2.3$~minutes. The additional selection binaries are largely determined once predictor confidences are computed, and contribute little to the branch-and-bound effort. Aggregating predictors does not aggregate their combinatorial cost, and the same conclusion applies to each embedding considered individually.

All four embeddings reduce the LP gap to values close to $2\%$, compared to over $6\%$ for the baseline formulation. This reflects the structure of the embeddings rather than their size: the deviation variables couple each decision variable to a learned target through a penalized $\ell_1$ term and act as a convex envelope, pulling the linear relaxation toward the learned feasibility pattern. The augmented formulations replace a weak polyhedral description of the operational feasible region with one informed by data. A tighter relaxation reduces the branch-and-bound effort, and offsets the cost introduced by the additional variables: the model grows without a proportional increase in solving time.
\section{Managerial insights}
\label{sec:managerial_insights}

\textbf{CLIF} is conceived as a decision-support system rather than as a replacement for human planners. In its intended deployment, the framework generates both the baseline plan $\mathbf{x}^\star$ produced by $\mathcal{LTRSP}_G$ and the augmented plan $\mathbf{x}^{\text{AUG}}$ produced by $\mathcal{LTRSP}_{\text{AUG}}^{\text{STK}}$, and presents both to the planner, who retains the final decision authority to choose and, if needed, adjust the plan. The role of the system is to provide a rule-compliant starting point that requires less manual modification than the baseline. The relevant performance metric is therefore whether \textbf{CLIF} plans satisfy the operational rules and reduce manual customization, both of which are verified in Section~\ref{sec:experiments}.

The penalty coefficient $\lambda$ controls the trade-off between routing and scheduling efficiency and adherence to the learned operational rules, and its calibration must be addressed before deployment. As shown in Figure~\ref{fig:all_rules}, the satisfaction rate reaches $100\%$ once $\lambda$ exceeds a rule-dependent threshold, after which the objective gap stabilizes. Daily retuning is not required. A single calibrated value of $\lambda^\star \approx 200$ achieves full satisfaction on all five rules while keeping the objective degradation below $1.7\%$. We recommend a one-time calibration on a held-out validation period, selecting the smallest value of $\lambda$ that achieves a target satisfaction rate of, say, $95\%$, and reusing this value in production. A coarse grid search over $\lambda \in \{50, 100, 200, 300, 500\}$ on fifty to one hundred historical instances has proven sufficient, and the calibration remains stable as long as the operational context does not change substantially.

The operational benefits of \textbf{CLIF} extend beyond satisfying hidden rules. Prior to deployment, planners spent on average between $1$ and $2$ hours per day manually adjusting baseline routes to comply with informal conventions. With \textbf{CLIF}, the augmented plan already incorporates the implicit rules learned from historical practice, so that the planner's role shifts from systematic rule enforcement to focused review and, when applicable, the introduction of new operational considerations not yet present in the data. Such manual edits are not lost: they are added to the training set and incorporated into subsequent plans through the next training cycle, so that the system progressively narrows the gap between the proposed plan and the executed one. Beyond the reduction in manual effort, this mechanism is expected to improve the consistency of executed plans across planners and across days, by enforcing a shared and explicitly maintained representation of operational conventions rather than relying on individual judgment. The combination of reduced manual workload and improved decision consistency constitutes the central practical benefit of the framework.

Integration into existing planning software is straightforward, as the framework is designed as a drop-in replacement for the baseline optimizer. The planning system first loads the daily demand, supply, and fleet data. The predictors $\hat{f}^{\text{LR}}_\theta$, $\hat{f}^{\text{DT}}_\theta$, and $\hat{f}^{\text{GNN}}_\theta$, trained offline on historical pairs $(\mathbf{x}_i^\star, \hat{\mathbf{x}}_i)$, are queried over all candidate keys to produce feasibility predictions. The stacked formulation $\mathcal{LTRSP}_{\text{AUG}}^{\text{STK}}$ is solved by Gurobi, and the resulting plan is presented to the planner for final review. The additional computation introduced by the embedding takes some additional minutes per instance Section~\ref{sec:computational_analysis}. Predictors are retrained weekly on the most recent months of executed plans, so that the learned feasibility patterns track 
evolving operational practice without intervention from the planning team. An online-learning extension, updating the predictors incrementally after each executed plan rather than through periodic batch retraining, would further reduce adaptation latency and is a natural direction for future work.

\section{Conclusion}
\label{sec:conclusion}
We presented \textbf{CLIF}, a framework that learns implicit operational constraints from historical routing executions and embeds the trained predictors exactly into a MILP formulation via linear constraints. Building on the methodological lineage of \cite{hewitt2020data} and \cite{bayani2024learning}, but departing from their reliance on a single model family, \textbf{CLIF} trains three complementary predictors, a GNN, decision trees, and linear regression, and combines them through a stacking mechanism within a single augmented optimization problem, $\mathcal{LTRSP}_{\mathrm{AUG}}$, that jointly enforces classical routing and scheduling constraints and learned feasibility structures without sacrificing solver guarantees.
On a large-scale industrial log-truck routing dataset, all four variants achieve $100\%$ satisfaction of the tested operational rules, so that the models are distinguished not by compliance but by the objective degradation they incur. Compared to the linear regression baseline of \cite{hewitt2020data} and the decision tree baseline of \cite{bayani2024learning}, the GNN embedding already achieves a lower objective gap on most rules by exploiting the relational structure of the routing and scheduling graph. The stacked model $\mathcal{LTRSP}_{\mathrm{AUG}}^{\mathrm{STK}}$ consistently achieves the lowest objective gap across all rules and under the combined multi-rule regime, confirming that combining predictors of different natures yields a closer approximation of the true feasible region. Since the nature of implicit rules is not known in advance in real-world operations, the stacking approach offers a robust choice for practical deployment.

Several limitations and directions for future work merit discussion. Validation against logged manual modifications, when these become available, remains a natural consolidation step. A pilot deployment of \textbf{CLIF} as a decision-support system with the industrial partner is currently under preparation and will provide complementary evidence on practitioner adoption and on the robustness of the framework to rules beyond the five categories considered here. Learning the penalty coefficient $\lambda$ jointly with the predictors and automating hyperparameter tuning across the three model families are also natural extensions. Finally, extending \textbf{CLIF} to an online learning setting, exploring additional predictor families within the stacking mechanism, and evaluating performance under larger and more diverse constraint regimes remain natural directions for future work.

\clearpage 
\begin{APPENDICES}

\section{Mathematical model of $\mathcal{LTRSP}$}
\label{LTRSP}
The $\mathcal{LTRSP}$ considers a heterogeneous fleet of trucks $V$, each assigned to a home base $h \in HB$ from which it departs and returns at the end of each operational day. Trucks transport wood products $p \in P$ from forest sites $f \in F$ to mills $m \in M$, subject to mill demands $d_{mp}$ and forest site supplies $s_{fp}$. Each truck performs a sequence of loaded trips (forest site to mill) and unloaded trips (home base to forest site, mill to forest site, or mill to home base) over a planning horizon $T$. The operating day is discretized into fixed time intervals $I$ to capture synchronization constraints arising from limited loader availability at forest sites and mills. This yields a space-time graph $\mathcal{G} = (\mathcal{N}, A)$, where nodes $\mathcal{N}$ are location-time pairs and arcs $A$ encode feasible movements, loading, unloading, and waiting operations. Arcs inconsistent with travel times or time windows are pruned prior to solving, reducing the problem size. The objective minimizes total hauling costs, waiting costs, and penalties for unmet demand, subject to flow conservation, loader capacity, time-window, supply-demand balance, and maximum trip constraints. The resulting Mixed-Integer Linear Program is detailed below.

\subsection{Mathematical notation}
We first provide a detailed overview of the mathematical notation used throughout this paper. The sets, indices, parameters and decision variables are introduced in Tables~\ref{tab1} and~\ref{deci-var}.
\begin{table}[h]
\tiny
\centering
\caption{Decision variables}
\begin{tabular}{|c|p{11.5cm}|}
\hline
\textbf{Decision variable} & \textbf{Description} \\
\hline
$x_{avt}$     & Binary variable equal to one if arc $a$ is used by truck $v$ on day $t$ \\
$\delta_{mp}$ & Unmet demand of product $p$ at mill $m$ \\
\hline
\end{tabular}
\label{deci-var}
\end{table}
\begin{table}[H]
\tiny
\centering
\caption{Sets and parameters}
\begin{tabular}{|c|p{11.5cm}|}
\hline  
\textbf{Notation} & \textbf{Description} \\
\hline
$F$               & Set of forest sites \\
$M$               & Set of mills \\
$V$               & Set of trucks \\
$P$               & Set of wood products \\
$I$               & Set of time intervals \\
$\mathcal{N}$     & Set of nodes at mills, forest sites, and home bases \\
$\mathcal{HB}_v$  & Set of possible home bases for truck $v$ \\
$T$               & Set of days \\
$A$               & Set of arcs, defined between two nodes at two time intervals \\
$A^+(n)$          & Set of outgoing arcs from node $n$ \\
$A^-(n)$          & Set of incoming arcs into node $n$ \\
$A_{f}^{mp}$      & Set of loaded arcs from forest site $f$ to mill $m$ transporting wood product $p$ \\
$\text{Source}_v$ & Start node for truck $v$ (departure from home base) \\
$\text{Sink}_v$   & End node for truck $v$ (return to home base) \\
$A_{fi}^{L}$      & Loading arcs at forest site $f$ at time interval $i$ \\
$A_{mi}^{U}$      & Unloading arcs at mill $m$ at time interval $i$ \\
$c_{a}^{h}$       & Hauling cost associated with arc $a$ \\
$c_{a}^{w}$       & Waiting cost associated with arc $a$ \\
$c_{mp}$          & Penalty cost of unmet demand of product $p$ at mill $m$ \\
$\tau_{a}^{+}$    & End interval time for arc $a$ \\
$\tau_{a}^{-}$    & Start interval time for arc $a$ \\
$d_{mp}$          & Demand of product $p$ at mill $m$ expressed in full truckloads \\
$q_{v}$           & Capacity of truck $v$ expressed in full truckloads \\
$K_{v}$           & Number of allowed trips for truck $v$ \\
$\beta_{v}$       & Parameter equal to $1$ if truck $v$ is not self-loader \\
$s_{fp}$          & Supply of product $p$ at forest site $f$ expressed in full truckloads \\
$S_{t}$           & Service time at forest sites (loading) and mills (unloading) \\
$\Lambda_{mi}$    & Number of available loaders at mill $m$ at time interval $i$ \\
$\Lambda_{fi}$    & Number of available loaders at forest site $f$ at time interval $i$ \\
$T_{a}^{v}$       & Travel time for arc $a$ \\
$T^o_{n},\, T^c_{n}$ & Opening and closing time at node $n$ \\
$\textit{BigM}$   & Very large constant \\
\hline
\end{tabular}
\label{tab1}
\end{table}

\subsection{Formulation}
The objective function and the constraints are defined over the sets, parameters, and decision variables introduced in Tables~\ref{tab1} and~\ref{deci-var}.
{\small
\begin{align}
\min \quad & \sum_{t,v,a} (c_{a}^{h}+c_{a}^{w}) x_{avt} 
             + \sum_{m,p} c_{mp} \delta_{mp} 
             \label{objfun} \\[5pt]
\text{s.t.} \quad 
& \sum_{a \in A^{+}(h)} x_{avt} \leq 1, 
  \quad \forall v, h \in \mathcal{HB}_v, t 
  \label{const1} \\[5pt]
& \sum_{a \in A^{-}(h)} x_{avt} = \sum_{a \in A^{+}(h)} x_{avt}, 
  \quad \forall v, h \in \mathcal{HB}_v, t 
  \label{const2} \\[5pt]
& \sum_{a \in A^+(n)} x_{avt} = \sum_{a \in A^-(n)} x_{avt}, 
  \quad \forall v, n, t 
  \label{const3} \\[5pt]
& \sum_{t,f,v} \sum_{a \in A_{f}^{mp}} q_v x_{avt} + \delta_{mp} = d_{mp}, 
  \quad \forall m, p 
  \label{const4} \\[5pt]
& \sum_{t,m,v} \sum_{a \in A_{f}^{mp}} q_v x_{avt} \leq s_{fp}, 
  \quad \forall f, p 
  \label{const5} \\[5pt]
& \sum_{t,m,f,p} \sum_{a \in A_{f}^{mp}} x_{avt} \leq K_v, 
  \quad \forall v 
  \label{const6} \\[5pt]
& \sum_{t,v} \sum_{a \in A_{fi}^L} \beta_v x_{avt} \leq \Lambda_{fi}, 
  \quad \forall i, f, t 
  \label{const7} \\[5pt]
& \sum_{t,v} \sum_{a \in A_{mi}^U} \beta_v x_{avt} \leq \Lambda_{mi}, 
  \quad \forall i, m 
  \label{const8} \\[5pt]
& \tau_a^{+} - M(1 - x_{avt}) + S_t \leq T^c_{an}, 
  \quad \forall v, a, t 
  \label{const9} \\[5pt]
& T^o_{an} \leq \tau_a^{+} + M(1 - x_{avt}), 
  \quad \forall v, a, t 
  \label{const10} \\[5pt]
& x_{avt} \in \{0, 1\}, 
  \quad \forall v, a, t 
  \label{const11} \\[5pt]
& \delta_{mp} \in \{0, \ldots, d_{mp}\}, 
  \quad \forall m, p 
  \label{const14}
\end{align}
}
The objective function~\eqref{objfun} minimizes the total cost incurred by all trips, combining hauling costs, waiting costs, and penalties for unmet mill demand. Constraint~\eqref{const1} ensures that each truck selects at most one starting arc from its home base. Constraints~\eqref{const2} and~\eqref{const3} enforce flow conservation at home bases and all intermediate nodes respectively, maintaining the continuity of each truck's route. Constraint~\eqref{const4} guarantees that mill demand is fully satisfied, up to the slack variable $\delta_{mp}$, while Constraint~\eqref{const5} ensures that the available supply at each forest site is not exceeded. Constraint~\eqref{const6} limits the total number of loaded trips per truck. Constraints~\eqref{const7} and~\eqref{const8} enforce loader availability at forest sites and mills respectively, ensuring that the number of trucks being served simultaneously does not exceed the number of available loaders at each time interval. Constraints~\eqref{const9} and~\eqref{const10} impose time-window compliance at each node, ensuring that trucks arrive and depart within the opening and closing times. Finally, Constraints~\eqref{const11} and~\eqref{const14} define the domain of the decision variables. This model is denoted $\mathcal{LTRSP}$.

\section{Local branching for perturbing $\mathcal{LTRSP}$ optimal plans}
\label{localbranching}
To generate an executed plan, we assume that a set of implicit operational constraints, denoted by $\mathcal{R}$, is applied to the optimal plan obtained from the baseline model $\mathcal{P}_m$. These implicit rules are modeled through additional linear constraints of the form $\mathbf{B}\,\mathbf{x} + \mathbf{C}\,\mathbf{y} \le \mathbf{e}$, which capture operational restrictions that are not explicitly included in the original planning model. Let $\mathcal{F}_{\mathcal{P}_m}$ denote the feasible region of $\mathcal{P}_m$, and let $\hat{x}_{atv}$ be the optimal solution of $\mathcal{P}_m$. The executed plan is generated by solving a local-branching problem that enforces the implicit constraints while remaining close to the original optimal plan. Specifically, we solve:
{\small
\begin{align}
\min_{\mathbf{x} \in \mathcal{F}_{\mathcal{P}_m}} \quad
& \sum_{t \in \mathcal{T}} \sum_{v \in \mathcal{V}} \sum_{a \in A}
  c_a\, x_{atv}, \label{eq:PLB_obj} \\[3pt]
\text{s.t.} \quad
& \mathbf{B}\,\mathbf{x} + \mathbf{C}\,\mathbf{y} \le \mathbf{e}, \label{eq:PLB_rules} \\[3pt]
& y_{atv} \ge x_{atv} - \hat{x}_{atv}, \quad \forall\, a,\, t,\, v, \label{eq:PLB_abs1} \\[3pt]
& y_{atv} \ge \hat{x}_{atv} - x_{atv}, \quad \forall\, a,\, t,\, v, \label{eq:PLB_abs2} \\[3pt]
& \sum_{t \in \mathcal{T}} \sum_{v \in \mathcal{V}} \sum_{a \in A}
  y_{atv} \le k. \label{eq:PLB_k}
\end{align}
}
The auxiliary variables $y_{atv}$ measure deviations from the optimal solution $\hat{x}_{atv}$, and $k$ defines the size of the local-branching neighborhood. The parameter $k$ controls the maximum number of arc-time-vehicle decisions that may differ from the optimal plan, thereby ensuring that the executed plan remains close to the optimal solution while satisfying the implicit operational constraints. In practice, the value of $k$ is tuned to balance realism in execution deviations and proximity to optimality.

\section{Linear regression and decision trees AOP models}
\label{AOP-trees-lineair}
The augmented optimization problem (AOP) formulations presented below follow the framework of \cite{hewitt2020data}. The augmented optimization model using linear regression is:
{\small
\begin{align}
& \min \quad \mathbf{c}^\top \mathbf{x} + \lambda \,\Delta, \label{eq:obj_simple} \\
& \text{s.t.} \quad \mathbf{A}\,\mathbf{x} \ge \mathbf{b}, \label{eq:lin_core} \\
& \hat{x}_\alpha - \Bigl(\hat{\beta}^\alpha + \textstyle\sum_{\alpha' \in \mathcal{A}} \beta_{\alpha'}^\alpha x_{\alpha'}\Bigr) = 0, \quad \forall \alpha \in \mathcal{A}, \label{eq:y_def} \\
& \delta_\alpha \ge \hat{x}_\alpha - x_\alpha, \quad \forall \alpha \in \mathcal{A}, \label{eq:abs_up} \\
& \delta_\alpha \ge x_\alpha - \hat{x}_\alpha, \quad \forall \alpha \in \mathcal{A}, \label{eq:abs_low} \\
& \Delta = \textstyle\sum_{\alpha \in \mathcal{A}} \delta_\alpha, \label{eq:Delta_def} \\
& \mathbf{x} \in \mathcal{X} \subset \{0,1\}^{|\mathcal{A}|},\quad \hat{\mathbf{x}} \in [ 0,1]^{|\mathcal{A}|}, \label{eq:domains_simple} \\
& \delta_\alpha \ge 0,\; \forall \alpha \in \mathcal{A}, \quad \Delta \ge 0. \label{eq:nonneg}
\end{align}
}
The augmented optimization model using decision trees, following the embedding of \cite{bayani2024learning}, is, where $\tilde{x}_\alpha^l := \hat{\beta}^{\alpha l} + \sum_{\alpha' \in \mathcal{A}} \beta_{\alpha'}^{\alpha l} x_{\alpha'}$ denotes the leaf-$l$ regression at arc $\alpha$, constraints \eqref{eq:beta_up}--\eqref{eq:right_branch} hold for all relevant $\alpha$, $l$, $m$, and constraints \eqref{eq:one_hot}--\eqref{eq:delta_def} for all $\alpha \in \mathcal{A}$:
{\small
\begin{align}
& \min \quad \mathbf{c}^\top \mathbf{x} + \lambda \,\Delta, \label{eq:main_obj} \\
& \text{s.t.} \quad \mathbf{A}\,\mathbf{x} \ge \mathbf{b}, \label{eq:main_core} \\
& \hat{x}_\alpha - \tilde{x}_\alpha^l \le M(1 - z^{\alpha l}), \label{eq:beta_up} \\
& \hat{x}_\alpha - \tilde{x}_\alpha^l \ge -M(1 - z^{\alpha l}), \label{eq:beta_low} \\
& \textstyle\sum_{\alpha' \in \mathcal{A}} h_{\alpha'}^{\alpha m} x_{\alpha'} \le g^{\alpha m} + (1 - z^{\alpha l}),\; m\!\in\!L, \label{eq:left_branch} \\
& \textstyle\sum_{\alpha' \in \mathcal{A}} h_{\alpha'}^{\alpha m}(x_{\alpha'} \!-\! e_{\alpha'}) \ge g^{\alpha m} \!-\! (1 \!-\! \epsilon)(1 \!-\! z^{\alpha l}),\; m\!\in\!R, \label{eq:right_branch} \\
& \textstyle\sum_{l \in L} z^{\alpha l} = 1, \label{eq:one_hot} \\
& \delta_\alpha \ge \hat{x}_\alpha - x_\alpha, \label{eq:abs_1} \\
& \delta_\alpha \ge x_\alpha - \hat{x}_\alpha, \label{eq:abs_2} \\
& \Delta = \textstyle\sum_{\alpha \in \mathcal{A}} \delta_\alpha, \label{eq:delta_def} \\
& \mathbf{x} \in \mathcal{X} \subset \{0,1\}^{|\mathcal{A}|},\; \hat{\mathbf{x}} \in [ 0,1]^{|\mathcal{A}|},\; z^{\alpha l} \in \{0,1\},\; \delta_\alpha \ge 0,\; \Delta \ge 0. \label{eq:domains}
\end{align}
}
\section{Instances characteristics}
\label{appendix:datasets}
Examples of instances on which the experiments are conducted are denoted by $\mathcal{W} = \{W_1, \ldots, W_{20}\}$, where each instance corresponds to a weekly operational planning scenario. The number of arcs is used as a proxy for network density, that is, the number of potential routes connecting supply and demand nodes. Network density is a key factor influencing computational effort, as it increases symmetry and degeneracy, thereby directly impacting the branching process. In addition, the problem structure is influenced by several operational characteristics, including a heterogeneous truck fleet, multiple product types, and diverse client demands. Table~\ref{instances-characteristics} summarizes the main characteristics of all instances.

\begin{table}[ht]
\centering
\caption{Characteristics of the 20 weekly instances. 
\textbf{Instance}: label of the weekly instance; 
\textbf{Dates}: planning horizon of the instance; 
\textbf{Served Mills}: number of mills served; 
\textbf{Used Blocks}: number of forest blocks used; 
\textbf{Products}: number of product types; 
\textbf{Vehicles}: number of available vehicles; 
\textbf{Home bases}: number of truck home bases; 
\textbf{Nodes}: number of nodes in the network; 
\textbf{Arcs}: number of arcs in the network; 
\textbf{Total Demand (GMT)}: total demand expressed in GMT units; 
\textbf{Avg Distance}: average distance traveled per arc (km); 
\textbf{Max Distance}: maximum distance between any two points (km).}
\resizebox{\textwidth}{!}{%
\begin{tabular}{cccccccccccc}
\hline
\textbf{Instance} & \textbf{Dates} & \textbf{Served Mills} & \textbf{Used Blocks} & \textbf{Products} & \textbf{Vehicles} & \textbf{Home bases} & \textbf{Nodes} & \textbf{Arcs} & \textbf{Total Demand (GMT)} & \textbf{Avg Distance} & \textbf{Max Distance} \\ \hline
$W_1$  & 2025-01-10 -- 2025-01-14 & 11 & 25 & 3 & 70 & 27 & 1295 & 26372 & 6140  & 242.06 & 574.00 \\ \hline
$W_2$  & 2024-12-30 -- 2025-01-04 & 13 & 24 & 3 & 60 & 19 & 1303 & 26504 & 8820  & 272.00 & 629.27 \\ \hline
$W_3$ & 2024-12-24 -- 2024-12-29 & 11 & 19 & 3 & 36 & 19 & 1063 & 17189 & 3820  & 259.31 & 651.28 \\ \hline
$W_4$  & 2024-11-29 -- 2024-12-03 & 13 & 23 & 3 & 51 & 25 & 1267 & 26934 & 7180  & 245.39 & 680.30 \\ \hline
$W_{5}$ & 2024-11-24 -- 2024-11-28 & 15 & 23 & 3 & 53 & 17 & 1321 & 29265 & 10440 & 206.22 & 628.94 \\ \hline
$W_6$ & 2024-12-04 -- 2024-12-08 & 14 & 24 & 3 & 57 & 25 & 1333 & 32603 & 9600  & 231.08 & 660.30 \\ \hline
$W_7$  & 2024-12-09 -- 2024-12-13 & 17 & 30 & 3 & 66 & 27 & 1629 & 45237 & 14250 & 227.31 & 574.37 \\ \hline
$W_8$ & 2024-11-19 -- 2024-11-23 & 14 & 22 & 3 & 67 & 27 & 1267 & 32044 & 10800 & 283.70 & 657.06 \\ \hline
$W_{9}$ & 2024-11-14 -- 2024-11-18 & 17 & 23 & 3 & 68 & 30 & 1376 & 37120 & 10320 & 229.56 & 660.30 \\ \hline
$W_{10}$ & 2024-11-09 -- 2024-11-13 & 16 & 21 & 3 & 62 & 26 & 1280 & 31172 & 9630  & 225.54 & 657.06 \\ \hline
$W_{11}$ & 2024-10-30 -- 2024-11-03 & 18 & 46 & 4 & 59 & 27 & 1465 & 44237 & 12930 & 241.16 & 728.76 \\ \hline
$W_{12}$ & 2024-10-14 -- 2024-10-18 & 16 & 26 & 4 & 56 & 26 & 1464 & 36989 & 8140  & 264.96 & 739.14 \\ \hline
$W_{13}$ & 2024-10-10 -- 2024-10-14 & 13 & 37 & 4 & 62 & 27 & 1777 & 40795 & 11330 & 252.04 & 739.14 \\ \hline
$W_{14}$  & 2025-01-05 -- 2025-01-09 & 16 & 31 & 3 & 88 & 27 & 1637 & 43155 & 20130 & 244.92 & 574.00 \\ \hline
$W_{15}$ & 2024-10-25 -- 2024-10-29 & 17 & 37 & 4 & 72 & 30 & 1888 & 53103 & 12850 & 246.67 & 657.06 \\ \hline
$W_{16}$ & 2024-12-19 -- 2024-12-23 & 17 & 35 & 3 & 71 & 28 & 1808 & 50238 & 10540 & 268.47 & 651.28 \\ \hline
$W_{17}$  & 2024-12-14 -- 2024-12-18 & 17 & 38 & 3 & 64 & 27 & 1917 & 52736 & 11630 & 257.75 & 659.57 \\ \hline
$W_{18}$ & 2024-11-04 -- 2024-11-08 & 18 & 46 & 4 & 72 & 27 & 2233 & 86959 & 20240 & 220.54 & 702.36 \\ \hline
$W_{19}$ & 2024-10-20 -- 2024-10-24 & 17 & 40 & 4 & 78 & 30 & 1996 & 60449 & 23500 & 234.08 & 721.61 \\ \hline
$W_{20}$ & 2024-10-15 -- 2024-10-19 & 16 & 36 & 4 & 85 & 27 & 1825 & 50790 & 17170 & 227.70 & 739.14 \\ \hline
\end{tabular}%
}
\label{instances-characteristics}
\end{table}

\section{Proof of the structural containment property}
\label{app:containment_proof}
We prove the inequalities stated in Remark~\ref{rem:containment}. For each predictor $m \in \mathcal{M} = \{\text{LR}, \text{DT}, \text{GNN}\}$, let $\mathcal{P}^m$ denote the feasible set of $\mathcal{LTRSP}_{\text{AUG}}^{m}$, with variables $(\mathbf{x}, \hat{\mathbf{x}}^m, \boldsymbol{\Delta})$ and $x_{\alpha} \in \{0,1\}$, and let $\mathcal{P}^m_{\text{LP}}$ be its LP relaxation, obtained by relaxing $x_{\alpha} \in \{0,1\}$ to $x_{\alpha} \in [0,1]$. Let $Z^m_{\text{LP}}(\lambda)$ denote the LP relaxation value at penalty $\lambda$. Similarly, $\mathcal{P}^{\text{STK}}$ is the feasible set of the stacked formulation, which adds the selection variables $\mathbf{s}$ with $s^m_{\alpha} \in \{0,1\}$; $\mathcal{P}^{\text{STK}}_{\text{LP}}$ is its LP relaxation, obtained by relaxing both $x_{\alpha}$ and $s^m_{\alpha}$ to $[0,1]$; and $Z^{\text{STK}}_{\text{LP}}(\lambda)$ is the associated LP relaxation value.
\begin{proof}
We prove the LP inequality; the integer inequality follows by the identical construction with $x_{\alpha} \in \{0,1\}$ instead of $x_{\alpha} \in [0,1]$. Fix $\lambda \geq 0$ and let $m^\star \in \arg\min_{m \in \mathcal{M}} Z^m_{\text{LP}}(\lambda)$, with optimal LP solution $(\mathbf{x}^\star, \hat{\mathbf{x}}^{m^\star,\star}, \boldsymbol{\Delta}^\star) \in \mathcal{P}^{m^\star}_{\text{LP}}$. We construct a feasible point $(\bar{\mathbf{x}}, \bar{\mathbf{x}}^{\text{STK}}, \bar{\boldsymbol{\Delta}}, \bar{\mathbf{s}}, \bar{\boldsymbol{\gamma}}^{\max}) \in \mathcal{P}^{\text{STK}}_{\text{LP}}$ with the same objective value.
Set $\bar{x}_{\alpha} = x^\star_{\alpha}$ for all $\alpha \in \mathcal{A}$. The predictions $\hat{x}^m_{\alpha} = \hat{f}^m_\theta(\mathcal{G},X)_{\alpha}$ depend only on the fixed graph $\mathcal{G}$ and feature matrix $X$, and are therefore identical across both formulations. Set $\bar{s}^{m^\star}_{\alpha} = 1$ and $\bar{s}^m_{\alpha} = 0$ for all $m \neq m^\star$ and all $\alpha \in \mathcal{A}$. The selection constraints $\sum_m \bar{s}^m_{\alpha} = 1$ and $\bar{s}^m_{\alpha} \in [0,1]$ are satisfied. The linking constraints
\[
\hat{x}^m_{\alpha} - M_c(1-\bar{s}^m_{\alpha}) \;\leq\; \bar{x}^{\text{STK}}_{\alpha} \;\leq\; \hat{x}^m_{\alpha} + M_c(1-\bar{s}^m_{\alpha})
\]
are satisfied by setting $\bar{x}^{\text{STK}}_{\alpha} = \hat{x}^{m^\star}_{\alpha}$, since for $m = m^\star$ the constraint reduces to equality, and for $m \neq m^\star$ the constraints are slack by $M_c$. The confidence variables $\bar{\gamma}^m_{\alpha} = |\hat{x}^m_{\alpha} - 0.5|$ are fixed by the predictions, and we set $\bar{\gamma}^{\max}_{\alpha} = \max_{m \in \mathcal{M}} \bar{\gamma}^m_{\alpha}$. The selection constraints $\bar{\gamma}^m_{\alpha} \geq \bar{\gamma}^{\max}_{\alpha} - M_\gamma(1-\bar{s}^m_{\alpha})$ are trivially satisfied for $m = m^\star$ and slack for $m \neq m^\star$. Finally, set $\bar{\Delta}_{\alpha} = \Delta^\star_{\alpha} \geq |x^\star_{\alpha} - \hat{x}^{m^\star}_{\alpha}| = |\bar{x}_{\alpha} - \bar{x}^{\text{STK}}_{\alpha}|$, which satisfies the deviation constraints.
The constructed point is feasible for $\mathcal{P}^{\text{STK}}_{\text{LP}}$ and achieves objective $\sum_{\alpha} c_{\alpha} \bar{x}_{\alpha} + \lambda \sum_{\alpha} \bar{\Delta}_{\alpha} = Z^{m^\star}_{\text{LP}}(\lambda)$. Hence $Z^{\text{STK}}_{\text{LP}}(\lambda) \leq Z^{m^\star}_{\text{LP}}(\lambda) = \min_{m \in \mathcal{M}} Z^m_{\text{LP}}(\lambda)$. Applying the same construction with $x_{\alpha} \in \{0,1\}$ yields $Z^{\text{STK}}(\lambda) \leq \min_{m \in \mathcal{M}} Z^m(\lambda)$.
\end{proof}
\section{Confidence-based stacking with majority-Voting Layer}
\label{app:voting}

We present here the full MILP formulation of the majority-voting layer added on top of the confidence-based stacking of Section~\ref{sec:stacking}. The confidence-based mechanism is preserved exactly as defined in $\mathcal{LTRSP}_{\mathrm{AUG}}^{\mathrm{STK}}$: for each predictive key $\alpha \in \mathcal{A}$, the most confident predictor among $\mathcal{M} = \{\text{LR}, \text{DT}, \text{GNN}\}$ is selected through $s_\alpha^m$, and its prediction is propagated to $\hat{x}_\alpha^{\mathrm{STK}}$. The voting layer adds three binary votes, one per predictor, and triggers an override whenever the vote of the confidence-selected model disagrees with the majority of the three votes. The final stacked prediction $\hat{x}_\alpha^{\mathrm{STK}}$ follows the confidence-based output when there is no disagreement, and follows the majority vote otherwise. All voting variables are determined endogenously by the predictions themselves, so that the solver cannot select votes strategically.

For each predictive key $\alpha \in \mathcal{A}$ and each predictor $m \in \mathcal{M}$, the voting layer introduces four additional sets of variables on top of those of $\mathcal{LTRSP}_{\mathrm{AUG}}^{\mathrm{STK}}$. The binary vote $b_\alpha^m \in \{0, 1\}$ is equal to $1$ if predictor $m$ predicts above $0.5$ and $0$ otherwise. The majority vote $V_\alpha \in \{0, 1\}$ is equal to $1$ if at least two predictors vote $1$. The product variable $w_\alpha^m \in \{0, 1\}$ encodes $s_\alpha^m \cdot b_\alpha^m$, that is, the vote of the confidence-selected model. The override indicator $o_\alpha \in \{0, 1\}$ is equal to $1$ if the confidence-selected vote differs from the majority vote.

The vote of each predictor is determined by its prediction through the indicator constraints
\begin{align}
b_\alpha^m = 1 &\implies \hat{x}_\alpha^m \geq 0.5, \quad \forall\, \alpha,\, m, \label{eq:vote_on}\\
b_\alpha^m = 0 &\implies \hat{x}_\alpha^m \leq 0.5 - \varepsilon, \quad \forall\, \alpha,\, m, \label{eq:vote_off}
\end{align}
where $\varepsilon > 0$ is a small strict-inequality margin. These constraints ensure that the votes are deterministic functions of the predictions, not free decisions of the solver. The majority vote $V_\alpha$ is then encoded by
\begin{align}
\sum_{m \in \mathcal{M}} b_\alpha^m &\geq 2 V_\alpha, \quad \forall\, \alpha, \label{eq:maj_lb}\\
\sum_{m \in \mathcal{M}} b_\alpha^m &\leq 1 + 2 V_\alpha, \quad \forall\, \alpha, \label{eq:maj_ub}
\end{align}
which together enforce $V_\alpha = 1 \iff \sum_{m} b_\alpha^m \geq 2$.

The product $w_\alpha^m = s_\alpha^m \cdot b_\alpha^m$ is linearized through the McCormick envelope
\begin{align}
w_\alpha^m &\leq s_\alpha^m, \quad \forall\, \alpha,\, m, \label{eq:mc1}\\
w_\alpha^m &\leq b_\alpha^m, \quad \forall\, \alpha,\, m, \label{eq:mc2}\\
w_\alpha^m &\geq s_\alpha^m + b_\alpha^m - 1, \quad \forall\, \alpha,\, m. \label{eq:mc3}
\end{align}
Since exactly one selection variable $s_\alpha^m$ is equal to one, the sum $\sum_{m} w_\alpha^m$ retrieves the vote of the confidence-selected model. The override indicator $o_\alpha$ is then equal to one if and only if this sum differs from the majority vote $V_\alpha$, which is enforced by
\begin{align}
o_\alpha &\geq V_\alpha - \sum_{m \in \mathcal{M}} w_\alpha^m, \quad \forall\, \alpha, \label{eq:ovr1}\\
o_\alpha &\geq \sum_{m \in \mathcal{M}} w_\alpha^m - V_\alpha, \quad \forall\, \alpha, \label{eq:ovr2}\\
o_\alpha &\leq V_\alpha + \sum_{m \in \mathcal{M}} w_\alpha^m, \quad \forall\, \alpha, \label{eq:ovr3}\\
o_\alpha &\leq 2 - V_\alpha - \sum_{m \in \mathcal{M}} w_\alpha^m, \quad \forall\, \alpha. \label{eq:ovr4}
\end{align}

The final stacked prediction $\hat{x}_\alpha^{\mathrm{STK}}$ follows the confidence-selected prediction when $o_\alpha = 0$, and follows the majority vote $V_\alpha$ when $o_\alpha = 1$. With $M_c = 1$ since all predictions lie in $[0, 1]$, the original coupling between $\hat{x}_\alpha^{\mathrm{STK}}$ and $\hat{x}_\alpha^m$ in $\mathcal{LTRSP}_{\mathrm{AUG}}^{\mathrm{STK}}$ is replaced by the hybrid coupling
\begin{align}
\hat{x}_\alpha^{\mathrm{STK}}
&\leq \hat{x}_\alpha^m + M_c(1 - s_\alpha^m) + M_c\, o_\alpha,
\quad \forall\, \alpha,\, m, \label{eq:hyb_conf_ub}\\
\hat{x}_\alpha^{\mathrm{STK}}
&\geq \hat{x}_\alpha^m - M_c(1 - s_\alpha^m) - M_c\, o_\alpha,
\quad \forall\, \alpha,\, m, \label{eq:hyb_conf_lb}\\
\hat{x}_\alpha^{\mathrm{STK}}
&\leq V_\alpha + M_c(1 - o_\alpha),
\quad \forall\, \alpha, \label{eq:hyb_vote_ub}\\
\hat{x}_\alpha^{\mathrm{STK}}
&\geq V_\alpha - M_c(1 - o_\alpha),
\quad \forall\, \alpha. \label{eq:hyb_vote_lb}
\end{align}
Constraints~\eqref{eq:hyb_conf_lb}--\eqref{eq:hyb_conf_ub} reduce to $\hat{x}_\alpha^{\mathrm{STK}} = \hat{x}_\alpha^m$ for the selected model $m$ when $o_\alpha = 0$, and become inactive when $o_\alpha = 1$. Conversely, constraints~\eqref{eq:hyb_vote_ub}-- \eqref{eq:hyb_vote_lb} reduce to $\hat{x}_\alpha^{\mathrm{STK}} = V_\alpha$ when $o_\alpha = 1$, and become inactive when $o_\alpha = 0$. The deviation $\Delta_\alpha = |x_\alpha - \hat{x}_\alpha^{\mathrm{STK}}|$ remains defined through the same linearization as in $\mathcal{LTRSP}_{\mathrm{AUG}}^{\mathrm{STK}}$ and enters the objective through the term $\lambda \sum_{\alpha} \Delta_\alpha$.
\end{APPENDICES}
\end{document}